\documentclass[12pt]{article} 
\author{Constantin N. Beli \footnote{Partially supported by the Contract
    2-CEx06-11-20.}}
\title{A new approach to classification of integral quadratic forms over
dyadic local fields\footnote{In [B1] this paper was announced under the title
``BONG version of O'Meara's 93:28  theorem. We changed the title at the
referee's suggestion.}}
\usepackage{amsfonts,amsmath,amssymb}

\def\a{\alpha} \def\b{\beta}   \def\D{\Delta} 
 \def\e{\varepsilon}  
 \def\h{\frac}  \def\j{\infty}
   \def\m{\lim}
   
\def\p{\partial}     
    \def\te{\theta}
   
\def\z{\longrightarrow} \def\({\overline}
\def\){\underline} \def\<{\cdot} \def\go{\mathfrak}
\def\>{~~~~~~~} \def\#{{\bf
Definition}} \def\*{\section} \def\be{\begin{equation}}
\def\ee{\end{equation}}

\def\sb{\subset} \def\sp{\supset} \def\sbq{\subseteq} \def\spq{\supseteq} 
\def\ti{\times} \def\od{{\rm ord}\,} \def\oo{{\cal O}} \def\pp{\perp}
\def\ss{{\go s}} \def\nn{{\go n}} \def\ww{{\go w}} \def\FF{{\go f}}
\def\GG{{\go g}} \def\ff{\dot{F}} \def\ooo{{\oo^\ti}} \def\AA{{\mathsf a}}
\def\BB{{\mathsf b}} \def\mo{{\rm mod}~}  \def\hh{{\rm H}}
\def\aa{A(0,0)} \def\ab{A(2,2\rho )} \def\fs{\ff^2}  
\def\p{\go p} \def\*{\sharp}  \def\0{}
 \def\1{^{-1}}  \def\dd{{\mathfrak d}}
\def\aaa{{\cal A}} \def\[{\prec} \def\]{\succ} 
\def\la{\langle} \def\ra{\rangle} 
\def\bmat{\left(\begin{array}} \def\emat{\end{array}\right)} \def\ev{\equiv}
\def\ap{\cong}   
\def\N{{\rm N}}

 \def\rep{{\rightarrow\!\!\!\! -}}
 \def\m2{~(\mo 2)} \def\no{\noindent}
 \def\bth{\begin{thm}} \def\eth{\end{thm}}
 \def\blem{\begin{lem}}
\def\elem{\end{lem}}

\newtheorem{theorem}{Theorem}[section]

\newtheorem{lemma}[theorem]{Lemma}
\newtheorem{definition}{Definition}
\newtheorem{corollary}[theorem]{Corollary}

\newtheorem{bof}[theorem]{}
\newtheorem{teorema}{Theorem}

\def\qed{\mbox{$\Box$}\vspace{\baselineskip}}
\def\pf{$Proof.$} 
\def\bco{\begin{corollary}} \def\eco{\end{corollary}} 
\def\bdf{\begin{definition}} \def\edf{\end{definition}} 
\def\btm{\begin{theorem}} \def\etm{\end{theorem}} 
\def\blm{\begin{lemma}} \def\elm{\end{lemma}} 
\def\bff{\begin{bof}\rm} \def\eff{\end{bof}}
\def\btr{\begin{teorema}} \def\etr{\end{teorema}}

\def\de{\newcommand} \de\tm[1]{{\no\bf Theorem~#1}} 
 \def\mb{\mathbb} 
\def\RR{{\mb R}}\def\QQ{{\mb Q}}\def\ZZ{{\mb Z}}\def\NN{{\mb
N}} 
\de\lm[1]{{\no\bf Lemma~#1}}
\de\df[1]{{\no\bf Definition~#1}} \de\co[1]{{\no\bf Corollary~#1}}
\de\tp[1]{\te (#1 )} \de\ts[1]{\te (O^-(#1 ))} \de\ty[1]{\te
(O(#1 ))} \de\tx[1]{\te (#1 )} \de\up[1]{(1+\p^{#1} )\fs}
 \de\upn[2]{(1+\p^{#1})\fs\cap\N (#2 )} \de\xt[2]{\te (#1 /#2 )}
 \de\ups[1]{((1+\p^{#1})\fs )^*} \de\upo[1]{(1+\p^{#1} )\ooo^2}
\de\upon[2]{(1+\p^{#1})\ooo^2\cap\N (#2 )}
\de\lr[1]{\longrightarrow^{\!\!\!\!\!\!\!\! #1}}
\de\lf[1]{\longleftarrow^{\!\!\!\!\!\!\!\! #1}}
\de\si[1]{\sim^{\!\!\!\!\! #1}} \de\apr[1]{\approx^{\!\!\!\!\! #1}}

\begin{document}
\maketitle
\begin{quote}
{\bf\footnotesize In [OM] O'Meara solved the classification problem for
lattices over dyadic local  fields in terms of Jordan
decompositions. In this paper we translate his result in terms of good
BONGs. BONGs (bases of norm generators) were introduced in [B] as a new
way of describing lattices over dyadic local fields. This result and
the notions we introduce here are a first step towards a solution of
the more difficult problem of representations of lattices over dyadic
fields.} 
\end{quote}

\section{Introduction}

Since the main result of this paper is given in terms of BONGs, which
were introduced in [B], we now give a reminder of some of the
definitions and results in that paper which we will use here. 

Throughout this paper $F$ is a dyadic local field, $\oo$ the ring of
integers, $\p$ the prime ideal, $\ooo :=\oo\setminus\p$ the group of
units, $e:=\od 2$ and $\pi$ is a fixed prime element. For $a\in\ff$ we
denote its quadratic defect by $\dd (a)$ and let $\D =1-4\rho$ be a
fixed unit with $\dd (\D)=4\oo$. 

We denote by $d:\ff/\fs\z\NN\cup\{\j\}$ the order of the ``{\em relative}
quadratic defect'' $d(a)=\od a\1\dd (a)$. If $a=\pi^R\e$, with $\e\in\ooo$,
then $d(a)=0$ if $R$ is odd and $d(a)=d(\e )=\od\dd (\e )$ if $R$ is
even. Thus $d(\ff )=\{ 0,1,3,\ldots,2e-1,2e,\j\}$. This function satisfies the
domination principle $d(ab)\geq\min\{ d(a),d(b)\}$. 

If $\a$ is a positive integer then $\up\a =\{ a\in\ff \mid d(a)\geq\a\}$ and
$\upo\a =\{ a\in\ooo \mid d(a)\geq\a\}$. For convenience we set $\up\a :=\{
a\in\ff\mid d(a)\geq\a\}$ and $\upo\a :=\{ a\in\ooo\mid d(a)\geq\a\}$ for any
$\a\in\RR\cup\{\j\}$. Thus $\up\a =\fs$ for $\a>2e$ and $\up\a =\ff$ for
$\a\leq 0$. If $d$ is the smallest element in $d(\ff )$ s.t. $\a\leq
d$ then $\up\a =\up d$. 

We denote by $(\cdot,\cdot )_\p :\ff/\fs\times\ff/\fs\z\{\pm 1\}$ the
Hilbert symbol, which is a non-degenerate bilinear symmetric form. 

If $a\in\ff$, we denote by $\N (a)$ the norm group $\N (F(\sqrt
a)/F)=\{ b\in\ff\mid (a,b)_\p =1\}$. If $b\in\ff$ and $d(a)+d(b)>2e$
then $(a,b)_\p =1$. However if $\a\notin\fs$ then there is $b\in\ff$
with $d(b)=2e-d(a)$ s.t. $(a,b)_\p =-1$. (For $d(a)$ odd this is just
[H, Lemma 3]. If $d(a)=2e$ and $b\in\ff$ is arbitrary with $d(b)=0$
then $a\in\D\fs$ and $\od b$ is odd so $(a,b)_\p =-1$. Similarly if
$d(a)=0$ and $d(b)=2e$ we have $(a,b)_\p =-1$.) Thus $\up\a\sbq\N (a)$
iff $\a +d(a)>2e$. 

An element $x$ of a lattice $L$ is called a {\em norm generator} of
$L$ if $\nn L=Q(x)\oo$. A sequence $x_1,\ldots,x_n$ of vectors in $FL$ is
called a {\em basis of norm generators} (BONG) for $L$ if $x_1$ is a
norm generator for $L$ and $x_2,\ldots,x_n$ is a BONG for
$pr_{x_1^\pp}L$. A BONG uniquely determines a lattice so, if
$x_1,\ldots,x_n$ is a BONG for $L$, we will write $L=\[
x_1,\ldots,x_n\]$. If moreover $Q(x_i)=a_i$ we say that $L\ap\[
a_1,\ldots,a_n\]$ relative to the BONG $x_1,\ldots,x_n$. If $L\ap\[
a_1,\ldots,a_n\]$ then $\det L=a_1\cdots a_n$. 

If $x_1,...,x_n$ are mutually orthogonal vectors with $Q(x_i)=a_i$, $L=\oo
x_1\pp\cdots\pp\oo x_n$ and $V=Fx_1\pp\cdots\pp Fx_n$ then we sat that
$L\ap\la a_1,\ldots,a_n\ra$ and $V\ap [a_1,\ldots,a_n]$ relative to the basis
$x_1,\ldots,x_n$. 

If $L$ is binary with $\nn L=\a\oo$, we denote by $a(L):=\det
L\,\a^{-2}$ and by $R(L):=\od vol L-2\od\nn L=\od
a(L)$. $a(L)\in\ff/\ooo^2$ is an invariant of $L$ and it determines the
class of $L$ up to scaling. If $L\ap\[\a,\b\]$ then $a(L)=\h\b\a$. 

We denote by $\aaa =\aaa_F\sb\ff/\ooo^2$ the set of all possible
values of $a(L)$, where $L$ is an arbitrary binary lattice. We have
$\aaa =\{ a\in\h 14\oo\mid a\neq 0,\dd (-a)\sbq\oo\}$. If $\od a=R$
and $d(-a)=d$, then $a\in\h 14\oo$ means $R\geq -2e$, while $\dd
(-a)\sbq\oo$ means $R+d=\od\dd (-a)\geq 0$. 

If $a(L)=a=\pi^R\e$ with $d(a)=d$ then:

$L$ is nonmodular, proper modular or improper modular iff $R>0$,
$R=0$, resp. $R<0$. 

If $R$ is odd then $R>0$. 

The inequality $R+2e\geq 0$ becomes equality iff $a\in -\h 14\ooo^2$ or
$a\in -\h\D 4\ooo^2$. We have $a(L)=-\h 14$ resp. $a(L)=-\h\D 4$ when
$L\ap\pi^r\aa$ resp. $\pi^r\ab$ for some integer $r$. 

The inequality $R+d\geq 0$ becomes equality iff $a\in -\h\D 4\ooo^2$. 

A special type of BONGs is the so called ``good BONGs''. If $L\ap\[
a_1,\ldots,a_n\]$ relative to some BONG $x_1,\ldots,x_n$ and $\od a_i=R_i$
we say that the BONG $x_1,\ldots,x_n$ is good if $R_i\leq R_{i+2}$ for
any $1\leq i\leq n-2$. 

{\bf Remark} The condition $R_i\leq R_{i+2}$ for $1\leq i\leq n-2$ is
equivalent to the condition that the sequence $(R_i+R_{i+1})$ is
increasing. 

A set $x_1,\ldots,x_n$ of orthogonal vectors with $Q(x_i)=a_i$ and $\od
a_i=R_i$ is a good BONG for some lattice iff $R_i\leq R_{i+2}$ for all
$1\leq i\leq n-2$ and $a_{i+1}/a_i\in\aaa$ for all $1\leq i\leq
n-1$. The condition $a_{i+1}/a_i\in\aaa$ is equivalent to
$R_{i+1}-R_i+2e\geq 0$ and $R_{i+1}-R_i+d(-a_ia_{i+1})\geq 0$. As
consequences of $a_{i+1}/a_i\in\aaa$, if $R_{i+1}-R_i$ is odd then it
is positive, if $R_{i+1}-R_i=-2e$ then $a_{i+1}/a_i\in -\h 14\ooo^2$
or $-\h\D 4\ooo^2$ and if $R_{i+1}-R_i+d(-a_ia_{i+1})=0$ then
$a_{i+1}/a_i\in -\h\D 4\ooo^2$. 

The good BONGs enjoy some properties similar to those of orthogonal
bases. If $L\ap\[ a_1,\ldots,a_n\]$ relative to some good BONG
$x_1,\ldots,x_n$ and $\od a_i=R_i$ then $L^\*\ap\[ a_1\1,\ldots,a_n\1\]$
relative to the good BONG $x_n^\*,\ldots,x_1^\*$, where $x_i^\* =Q(x)\1
x_i$. Also if for some $1\leq i\leq j\leq n$ we have $\[
x_i,\ldots,x_j\]\ap\[ b_i,\ldots,b_j\]$ relative to some other good BONG
$y_i,\ldots,y_j$ then $L\ap\[
a_1,\ldots,a_{i-1},b_i,\ldots,b_j,a_{i+1},\ldots,a_n\]$ relative to the good
BONG $x_1,\ldots,x_{i-1},y_i,\ldots,y_j,x_{i+1},\ldots,x_n$. There are some
differences though from the orthogonal bases. E.g. the relation $L=\[
x_1,\ldots,x_i\]\pp\[ x_{i+1},\ldots,x_n\]$ holds iff $R_i\leq R_{i+1}$. 

The orders $R_i=\od a_i$ are independent of the choice of the good
BONGs and they are in 1-1 correspondence with the invariants $t,\dim
L_k,\ss_k:=\ss L_k$ and $\nn L^{\ss_k}$, where $L=L_1\pp\ldots\pp L_t$ is
a Jordan splitting. More precisely, if $\ss_k=\p^{r_k}$, $\nn
L^{\ss_k}=\p^{u_k}$ and $n_k=\sum_{l\leq k}\dim L_l$, then the sequence
$R_{n_{k-1}+1},\ldots,R_{n_k}$ is $r_k,\ldots,r_k$ if $L_k$ is proper
(i.e. if $r_k=u_k$), and it is $u_k,2r_k-u_k,\ldots,u_k,2r_k-u_k$
otherwise; see [B, Lemma 4.7].

The good BONGs are closely connected with the {\em maximal norm
splittings}. A splitting $L=L_1\pp\ldots\pp L_t$ is called a maximal norm
splitting if $\ss L_1\spq\ldots\spq\ss L_t$ and $\dim L_i\leq 2$, $L_i$
is modular and $\nn L_i=\nn L^{\ss L_i}$ for all $1\leq i\leq
t$. Condition $\nn L_i=\nn L^{\ss L_i}$ is equivalent to $\nn
L_1\spq\ldots\spq\nn L_t$ and $\nn L_1^\*\sbq\ldots\sbq\nn L_t^\*$. If we
put together the BONGs of the components $L_1,\ldots,L_t$ of a maximal
norm splitting we get a good BONG for $L$. Conversely any good BONG of
a lattice can be obtained by putting together some BONGs of the
components of some maximal norm splitting. Moreover, the splitting can
be chosen s.t. all binary components are improper modular. An explicit
algorithm for finding a maximal norm splitting and, hence, a good BONG of a
lattice is provided in [B1, Section 7]. 

\section{The invariants $\a_i$}

Let $L$ be a lattice over the dyadic field $F$. Let $L\ap\[
a_1,\ldots,a_n\]$ relative to a good BONG and let $R_i:=\od a_i$. Also
let $L=L_1\pp\ldots\pp L_t$ be a Jordan decomposition. We keep the
notations of [OM] $\ss_k:=\ss L_k$, $\GG_k:=\GG L^{\ss_k}$, $\ww_k:=\ww
L^{\ss_k}$ but, in order to avoid confusion, we write $\AA_k$ for
O'Meara's $a_k$. Also we denote $r_k=\od\ss_k$, $u_k=\od\AA_k=\od\nn
L^{\ss_k}$, Associated to our splitting we have the Jordan chain
$L_{(1)}\sb\ldots\sb L_{(t)}$ and the inverse Jordan chain
$L_{(1)}^*\sp\ldots\sp L_{(t)}^*$, where $L_{(k)}:=L_1\pp\ldots\pp L_k$ and
$L_{(k)}^*:=L_k\pp\ldots\pp L_t$. 

Since $R_i$'s are invariants of $L$ we will write $R_i=R_i(L)$. 

\bdf For any $1\leq i\leq n-1$ we define $\a_i=\a_i(L)$ by: 
\begin{multline*}
\a_i:=\min (\{(R_{i+1}-R_i)/2+e\}\cup\{
R_{i+1}-R_j+d(-a_ja_{j+1})\mid 1\leq j\leq i\}\\ \cup\{
R_{j+1}-R_i+d(-a_ja_{j+1})\mid i\leq j<n\}). 
\end{multline*}
\edf

Apparently $\a_i(L)$ defined this way depends on the choice of the
good BONG. We will show later that, in fact, it depends only on
$L$. For the time being we will mean $\a_i(L)$ with respect to a given
good BONG. We give now some properties of $\a_i$'s.

\blm If $k\leq i<l$ then, in the set defining $\a_i$, we can replace
$(R_{i+1}-R_i)/2+e$ and all the terms corresponding to indices $k\leq j<l$, by
$\a_{i-k+1}(\[ a_k,\ldots,a_l\] )$. In particular, $\a_i\leq\a_{i-k+1}(\[
a_k,\ldots,a_l\] )$. 
\elm
\pf By definition $\a_{i-k+1}(\[ a_k,\ldots,a_l\] )=\min
(\{(R_{i+1}-R_i)/2+e\}\cup\{ R_{i+1}-R_j+d(-a_ja_{j+1})\mid k\leq j\leq
i\}\cup\{ R_{j+1}-R_i+d(-a_ja_{j+1})\mid i\leq j< l\})$. Hence the
conclusion. \qed

\blm The sequence $(R_i+\a_i)$ is increasing and the sequence
$(-R_{i+1}+\a_i)$ is decreasing.
\elm

\pf Let $1\leq i\leq h\leq n-1$. We have $R_i+R_{i+1}\leq
R_h+R_{h+1}$. From Definition 1 we get $R_i+\a_i=\min (\{
(R_i+R_{i+1})/2+e\}\cup\{ R_i+R_{i+1}-R_j+d(-a_ja_{j+1})\mid 1\leq j\leq
i\}\cup\{ R_{j+1}+d(-a_ja_{j+1})\mid i\leq j<n\} )$ and
$-R_{i+1}+\a_i=\min (\{ -(R_i+R_{i+1})/2+e\}\cup\{
-R_j+d(-a_ja_{j+1})\mid 1\leq j\leq i\}\cup\{
R_{j+1}-R_i-R_{i+1}+d(-a_ja_{j+1})\mid i\leq j<n\} )$. Similarly for
$R_h+\a_h$ and $-R_{h+1}+\a_h$. In order to prove that $R_i+\a_i\leq
R_h+\a_h$ we show that the elements in the set that has $R_i+\a_i$ as
its minimum are $\leq$ than the corresponding elements for
$R_h+\a_h$. Same for $-R_{i+1}+\a_i\geq -R_{h+1}+\a_h$. 

The proof is straightforward and uses the fact that $R_l+R_{l+1}$ is an
increasing sequence. For terms involving $d(-a_ja_{j+1})$ we consider the
cases $j\leq i$, $i\leq j\leq h$ and $h\leq j$ and use the inequalities among
$R_i+R_{i+1}$, $R_j+R_{j+1}$ and $R_h+R_{h+1}$ that occur in each case.\qed

\bco Suppose that $1\leq i\leq j\leq n-1$ and
$R_i+R_{i+1}=R_j+R_{j+1}$. Then: 

(i) $R_i+\a_i=\ldots=R_j+\a_j$ and $-R_{i+1}+\a_i=\ldots=-R_{j+1}+\a_j$. 

(ii) $R_k=R_l$ for any $k,l\in [i,j+1]$ of the same parity and
$\a_k=\a_l$ for any $k,l\in [i,j]$ of the same parity. 

(iii) If $\a_k=(R_{k+1}-R_k)/2+e$ for some $i\leq k\leq j$ then
$\a_k=(R_{k+1}-R_k)/2+e$ for all $i\leq k\leq j$. 

In the particular case when $j=i+1$ we get the following statement:

If $1\leq i\leq n-2$ and $R_i=R_{i+2}$ then $R_i+\a_i=R_{i+1}+\a_{i+1}$,
$-R_{i+1}+\a_i=-R_{i+2}+\a_{i+1}$ and $\a_i=(R_{i+1}-R_i)/2+e$ is equivalent
to $\a_{i+1}=(R_{i+2}-R_{i+1})/2+e$. 
\eco
\pf For (i) we note that $R_i+R_{i+1}=(R_i+\a_i)-(-R_{i+1}+\a_i)$ and
$R_j+R_{j+1}=(R_j+\a_j)-(-R_{j+1}+\a_j)$ and use Lemma 2.2. By using the
fact that $R_k+R_{k+1}$ is an increasing sequence we get
$R_i+R_{i+1}=R_{i+1}+R_{i+2}=\ldots =R_j+R_{j+1}$, which is equivalent to
(ii). Finally (iii) follows from $R_i+\a_i=\ldots =R_j+\a_j$,
$R_i+R_{i+1}=\ldots =R_j+R_{j+1}$ and the fact that $\a_k=(R_{k+1}-R_k)/2+e$ is
equivalent to $R_k+\a_k=(R_k+R_{k+1})/2+e$.\qed

\blm Suppose that $1\leq i<n$ and $1\leq k\leq h<l\leq n$. Then: 

(i) If $h\leq i$ then all terms in the definition of $\a_i$
corresponding to indices $k\leq j\leq h$ can be replaced by
$R_{i+1}-R_{h+1}+\a_{h-k+1}(\[ a_k,\ldots,a_l\] )$. In particular, all
terms with $1\leq j\leq h$ can be replaced by $R_{i+1}-R_{h+1}+\a_h$. 

(ii) If $i\leq h$ then all terms in the definition of $\a_i$
corresponding to indices $h\leq j<l$ can be replaced by
$R_h-R_i+\a_{h-k+1}(\[ a_k,\ldots,a_l\] )$. In particular, all terms with
$h\leq j<n$ can be replaced by $R_h-R_i+\a_h$. 
\elm
\pf By Lemma 2.1 we have $\a_{h-k+1}(\[ a_k,\ldots,a_l\] )\geq\a_h$. 

(i) By Lemma 2.2 we have $\a_i\leq R_{i+1}-R_{h+1}+\a_h\leq
R_{i+1}-R_{h+1}+\a_{h-k+1}(\[ a_k,\ldots,a_l\] )$. If $k\leq j\leq h$
then $\a_{h-k+1}(\[ a_k,\ldots,a_l\] )\leq R_{h+1}-R_j+d(-a_ja_{j+1})$ so
$R_{i+1}-R_{h+1}+\a_{h-k+1}(\[ a_k,\ldots,a_l\] )\leq
R_{i+1}-R_j+d(-a_ja_{j+1})$. Therefore if we add
$R_{i+1}-R_{h+1}+\a_{h-k+1}(\[ a_k,\ldots,a_l\] )$ to the set that
defines $\a_i$ and remove any one of $R_{i+1}-R_j+d(-a_ja_{j+1})$ with
$k\leq j\leq h$ then $\a_i$ does not change. 

(ii) By Lemma 2.2 we have $\a_i\leq R_h-R_i+\a_h\leq
R_h-R_i+\a_{h-k+1}(\[ a_k,\ldots,a_l\] )$. If $h\leq j<l$ then
$\a_{h-k+1}(\[ a_k,\ldots,a_l\] )\leq R_{j+1}-R_h+d(-a_ja_{j+1})$ so
$R_h-R_i+\a_{h-k+1}(\[ a_k,\ldots,a_l\] )\leq
R_{j+1}-R_i+d(-a_ja_{j+1})$. Thus if we add $R_h-R_i+\a_{h-k+1}(\[
a_k,\ldots,a_l\] )$ to the set that defines $\a_i$ and remove any one of
$R_{j+1}-R_i+d(-a_ja_{j+1})$ with $h\leq j<l$ then $\a_i$ does not change. 

If we take $k=1$ and $l=n$ then $\a_{h-k+1}(\[ a_k,\ldots,a_l\] )$
becomes $\a_h(\[ a_1,\ldots,a_n\] )=\a_h(L)=\a_h$ so we get the second
claims of (i) and (ii). \qed

\bco For any $1\leq i\leq n-1$ we have: 

(i) $\a_i=\min\{ (R_{i+1}-R_i)/2+e,R_{i+1}-R_i+d(-a_ia_{i+1}),
R_{i+1}-R_i+\a_{i-1},R_{i+1}-R_i+\a_{i+1}\}$. 

(ii) $\a_i=\min\{ (R_{i+1}-R_i)/2+e,R_{i+1}-R_i+d(-a_ia_{i+1}),
R_{i+1}-R_i+\a_{i-1}(\[a_1,\ldots,a_i\] ),R_{i+1}-R_i+\a_1(\[
a_{i+1},\ldots,a_n\] )\}$. 

(The terms that do not make sense, i.e. $R_{i+1}-R_i+\a_{i-1}$ and
$R_{i+1}-R_i+\a_{i-1}(\[a_1,\ldots,a_i\] )$ when $i=1$, or
$R_{i+1}-R_i+\a_{i+1}$ and $R_{i+1}-R_i+\a_1(\[ a_{i+1},\ldots,a_n\] )$ when
$i=n-1$, are ignored.) 
\eco

\pf (i) By Lemma 2.4 (i) resp. (ii), in the set defining $\a_i$,
$R_{i+1}-R_i+\a_{i-1}$ can replace all the terms
$R_{i+1}-R_j+d(-a_ja_{j+1})$ with $1\leq j\leq i-1$, while
$R_{i+1}-R_i+\a_{i+1}$ replaces all $R_{j+1}-R_i+d(-a_ja_{j+1})$ with
$i+1\leq j<n$. Therefore $\a_i=\min\{ (R_{i+1}-R_i)/2+e,
R_{i+1}-R_i+d(-a_ia_{i+1}),R_{i+1}-R_i+\a_{i-1},R_{i+1}-R_i+\a_{i+1}\}$.

(ii) Same as (i) but this time the terms corresponding to $1\leq j\leq
i-1$ are replaced by $R_{i+1}-R_i+\a_{i-1}(\[a_1,\ldots,a_i\] )$ and
those corresponding to $i+1\leq j<n$ by $R_{i+1}-R_i+\a_1(\[
a_{i+1},\ldots,a_n\] )$. \qed

\bff {\bf Remark} We have $L^\*\ap\[ a^\*_1,\ldots,a^\*_n\]$ with
$a^\*_i=a_{n+1-i}\1$ and $R^\*_i:=\od a^\*_i=-R_{n+1-i}$. One can easily see
that $\a^\*_i:=\a_i(L^\* )=\a_{n-i}$. Also $\a_i$'s are invariant to scaling. 
\eff

\blm If $1\leq i\leq n-1$ then:

(i) $\a_i\geq 0$ with equality iff $R_{i+1}-R_i=-2e$. 

(ii) If $R_{i+1}-R_i\geq 2e$ then $\a_i=(R_{i+1}-R_i)/2+e$. 

(iii) If $R_{i+1}-R_i\leq 2e$ then $\a_i\geq R_{i+1}-R_i$ with
equality iff $R_{i+1}-R_i=2e$ or it is odd. 

(iv) $\a_i$ is an odd integer unless $\a_i=(R_{i+1}-R_i)/2+e$. 
\elm

\pf We use induction on $n$. For $n=1$ our lemma is vacuous.

For the induction step let $1\leq i\leq n-1$ and let $L'=\[
a_1,\ldots,a_i\]$ and $L''=\[ a_{i+1},\ldots,a_n\]$. By Corollary 2.5(ii) we
have $\a_i=\min\{(R_{i+1}-R_i)/2+e,R_{i+1}-R_i+d(-a_1a_2),
R_{i+1}-R_i+\a,R_{i+1}-R_i+\b\}$, where $\a=\a_{i-1}(L')$ and
$\b=\a_1(L'')$. (We ignore $\a$ and $\b$ whenever they are not defined.) By
the induction hypothesis $\a,\b$ satisfy (i)-(iv) of the lemma. 

We have $(R_{i+1}-R_i)/2+e\geq 0$ with equality iff $R_{i+1}-R_i=-2e$
and $R_{i+1}-R_i+d(-a_ia_{i+1})\geq 0$ with equality iff $a_{i+1}/a_i\in
-\h\D 4\ooo^2$ which implies $R_{i+1}-R_i=-2e$. If
$R_{i+2}-R_{i+1}>2e$ then $\b =(R_{i+2}-R_{i+1})/2+e>2e$ so
$R_{i+1}-R_i+\b>R_{i+1}-R_i+2e\geq 0$. Similarly with $R_{i+1}-R_i+\a$
if $R_i-R_{i-1}>2e$. If $R_{i+2}-R_{i+1}\leq 2e$ then, by the
induction hypothesis, $\b\geq R_{i+2}-R_{i+1}$ with equality iff
$R_{i+2}-R_{i+1}$ is odd or it is $2e$. Thus $R_{i+1}-R_i+\b\geq
R_{i+2}-R_i\geq 0$ with equality iff $R_i=R_{i+2}$ and
$R_{i+2}-R_{i+1}$ is odd or $2e$. Suppose this happens. If
$R_{i+2}-R_{i+1}=2e$ then $R_{i+1}-R_i=R_{i+1}-R_{i+2}=-2e$. If
$R_{i+2}-R_{i+1}$ is odd then so is $R_{i+1}-R_i=R_{i+1}-R_{i+2}$ so
both must be positive. But this is impossible. Similarly for
$R_{i+1}-R_i+\a$ when $R_i-R_{i-1}\leq 2e$. Thus we have (i). 

If $R_{i+1}-R_i\geq 2e$ then $\a,\b\geq 0$ so
$R_{i+1}-R_i+d(-a_1a_2),R_{i+1}-R_i+\a,R_{i+1}-R_i+\b\geq
R_{i+1}-R_i\geq (R_{i+1}-R_i)/2+e$. Hence $\a_i=(R_{i+1}-R_i)/2+e$ and
we have (ii). 

We prove now (iii). If $R_{i+1}-R_i=2e$ then (ii) implies that
$\a_i=(R_{i+1}-R_i)/2+e=2e=R_{i+1}-R_i$ so we are done. If
$R_{i+1}-R_i<2e$ is odd then $d(-a_ia_{i+1})=0$ and $\a,\b\geq 0$ so
$\a_i=\min\{ (R_{i+1}-R_i)/2+e,R_{i+1}-R_i\} =R_{i+1}-R_i$. Finally if
$R_{i+1}-R_i<2e$ is even then $\od a_ia_{i+1}=R_i+R_{i+1}$ is even so
$d(-a_ia_{i+1})>0$. Also $R_i-R_{i-1},R_{i+2}-R_{i+1}\geq
R_i-R_{i+1}>-2e$ ($R_{i-1}\leq R_{i+1}$ and $R_i\leq R_{i+2}$) so by
(i) $\a,\b>0$. We have
$R_{i+1}-R_i+d(-a_1a_2),R_{i+1}-R_i+\a,R_{i+1}-R_i+\b
>R_{i+1}-R_i$. Since also $(R_{i+1}-R_i)/2+e>R_{i+1}-R_i$ (we have
$R_{i+1}-R_i<2e$) we get $\a_i>R_{i+1}-R_i$. 

We prove now (iv). If $R_{i+1}-R_i\geq 2e$ then (ii) implies
$\a_i=(R_{i+1}-R_i)/2+e$ so (iv) is vacuous. If $R_{i+1}-R_i<2e$ is odd
then (iii) implies $\a_i=R_{i+1}-R_i$ so $\a_i$ is odd. If
$R_{i+1}-R_i<2e$ is even then again $\od a_ia_{i+1}$ is even so
$d(-a_ia_{i+1})>0$. Suppose $\a_i<(R_{i+1}-R_i)/2+e$. If
$\a_i=R_{i+1}-R_i+d(-a_ia_{i+1})$ then if $d(-a_ia_{i+1})$ is odd
$\a_i$ will also be odd so we are done. Otherwise $d(-a_ia_{i+1})=2e$
or $\j$ so $\a_i=R_{i+1}-R_i+d(-a_ia_{i+1})\geq
R_{i+1}-R_i+2e\geq (R_{i+1}-R_i)/2+e>\a_i$. (We have
$R_{i+1}-R_i+2e\geq 0$.) Contradiction. If $\a_i=R_{i+1}-R_i+\a$ then
$\a_i$ is odd unless $\a$ is not odd which would imply
$\a=(R_i-R_{i-1})/2+e$. So $\a_i=R_{i+1}-R_i+(R_i-R_{i-1})/2+e\geq
(R_{i+1}-R_i)/2+e>\a_i$. (We have $R_{i+1}\geq
R_{i-1}$.) Contradiction. Similarly if $\a_i=R_{i+1}-R_i+\b$ since
$R_{i+1}-R_i+(R_{i+2}-R_{i+1})/2+e\geq (R_{i+1}-R_i)/2+e>\a_i$. (We
have $R_{i+2}\geq R_i$.) \qed

\bco (i) $\a_i\in\ZZ$ except when $R_{i+1}-R_i$ is odd and $>2e$. 

(ii) $\a_i$ is $<2e$, $=2e$ or $>2e$ if $R_{i+1}-R_i$ is $<2e$, $=2e$ or $>2e$
accordingly. 

(iii) $\a_i\in ([0,2e]\cap\ZZ )\cup ((2e,\j)\cap\h 12\ZZ )$. 
\eco
\pf (i) If $R_{i+1}-R_i>2e$ then $\a_i=(R_{i+1}-R_i)/2+e$. If
$R_{i+1}-R_i$ is even then $\a_i\in\ZZ$, while if it is odd then
$\a_i\in\h 12\ZZ\setminus\ZZ$. Suppose now that $R_{i+1}-R_i\leq
2e$. If $R_{i+1}-R_i$ is odd then $\a_i=R_{i+1}-R_i\in\ZZ$. If
$R_{i+1}-R_i$ is even then either $\a_i$ is an odd integer or
$\a_i=(R_{i+1}-R_i)/2+e\in\ZZ$. 

(ii) If $R_{i+1}-R_i<2e$ then $\a_i\leq (R_{i+1}-R_i)/2+e<2e$. If
$R_{i+1}-R_i=2e$ then $\a_i=(R_{i+1}-R_i)/2+e=2e$. If $R_{i+1}-R_i>2e$
then $\a_i=(R_{i+1}-R_i)/2+e>2e$. 

(iii) We have $\a_i\geq 0$. If $\a_i\leq 2e$ then $R_{i+1}-R_i\leq 2e$
so $\a_i\in\ZZ$. If $\a_i>2e$ then $R_{i+1}-R_i>2e$ so
$\a_i=(R_{i+1}-R_i)/2+e\in (2e,\j )\cap\h 12\ZZ$. \qed

\bco In each of the following cases, $\a_i$ depends only on
$R_{i+1}-R_i$: 

(i) If $R_{i+1}-R_i\geq 2e$ or $R_{i+1}-R_i\in\{ -2e,2-2e,2e-2\}$ then
$\a_i=(R_{i+1}-R_i)/2+e$. 

(ii) If $R_{i+1}-R_i$ is odd, then $\a_i=\min\{
(R_{i+1}-R_i)/2+e,R_{i+1}-R_i\}$. 
\eco
\pf (i) If $R_{i+1}-R_i\geq 2e$ then $\a_i=(R_{i+1}-R_i)/2+e$ by Lemma
2.7(ii). If $R_{i+1}-R_i=-2e$ then $\a_i=0=(R_{i+1}-R_i)/2+e$. If
$R_{i+1}-R_i=2-2e$ then $\a_i\in\ZZ$ and $0<\a_i\leq
(R_{i+1}-R_i)/2+e=1$ so $\a_i=1=(R_{i+1}-R_i)/2+e$. If
$R_{i+1}-R_i=2e-2$ then $\a_i\in\ZZ$ and $2e-2=R_{i+1}-R_i<\a_i\leq
(R_{i+1}-R_i)/2+e=2e-1$ so $\a_i=2e-1=(R_{i+1}-R_i)/2+e$. 

(ii) We use Lemma 2.7(ii) and (iii). If $R_{i+1}-R_i>2e$ then
$\a_i=(R_{i+1}-R_i)/2+e<R_{i+1}-R_i$. If $R_{i+1}-R_i<2e$ then
$\a_i=R_{i+1}-R_i<(R_{i+1}-R_i)/2+e$. In both cases $\a_i=\min\{
(R_{i+1}-R_i)/2+e,R_{i+1}-R_i\}$. \qed

\blm Let $\AA$ be a norm generator of a lattice $L$ and let $\ww\spq 2\ss L$
be a fractional ideal. Then $\ww =\ww L$ iff $\GG L=\AA\oo^2+\ww$ and we
have either $\ww =2\ss L$ or $\od\AA +\od\ww$ is odd. 
\elm
\pf For the necessity see [OM, 93A]. For the sufficiency it is enough to prove
that, given another fractional ideal $\ww'$ satisfying the hypothesis of the
lemma, we have $\ww =\ww'$. Suppose that $\ww\neq\ww'$. We may assume that
$\ww\sp\ww'$. Since $\ww\sp\ww'\spq 2\ss L$ we must have that $\od\AA +\od\ww$
is odd. Let $\ww =\BB\oo$. Then $\AA +\BB\in\AA\oo^2+\ww =\GG
L=\AA\oo^2+\ww'$. So $\AA +\BB =\AA\a^2+\BB'$ for some $\a\in\oo$ and
$\BB'\in\ww'\sb\ww$. It follows that $1+\BB/\AA =\a^2+\BB'/\AA$, which implies
that $\dd (1+\BB/\AA )\sbq\BB'/\AA\oo\sb\AA\1\ww$. On the other hand
$\od\BB/\AA =\od\AA\1\ww$ is odd and, since $\BB\oo =\ww\sbq\GG L\sbq\AA\oo$
and $\BB\oo =\ww\sp 2\ss L\spq 4\AA\oo$, we have
$4\oo\sb\BB/\AA\oo\sbq\oo$. By [OM, 63:5] we get $\dd (1+\BB/\AA )=\BB/\AA\oo
=\AA\1\ww$. Contradiction. \qed

\blm Let $J_1,\ldots,J_s$ be lattices in the same quadratic space and
let $J=\sum J_k$. If $\AA_k$ and $\AA$ are norm generators for $J_k$ and
$\AA$ and $J$, respectively, then: 
$$\GG J=\sum\GG J_k+2\ss J\text{ and }\ww J=\sum\ww
J_k+\sum\AA\1\dd (\AA\AA_k)+2\ss J.$$ 
\elm
\pf We have $\GG J_k\sbq\GG J$ and $2\ss J\sbq\GG J$ so $\GG
J\spq\sum\GG J_k+2\ss J$. For the reverse inclusion note that
$Q(J)\sbq\sum Q(J_k)+2\ss J$. Thus $\GG J=Q(J)+2\ss J\sbq\sum
(Q(J_k)+2\ss J_k)+2\ss J=\sum\GG J_k+2\ss J$. 

We have $\AA\oo^2\sbq\GG J$ and $2\AA\oo =2\nn J\sbq 2\ss J\sbq\GG J$
so $\GG J=\AA\oo^2+2\AA\oo +\GG J=\AA\oo^2+2\AA\oo +\sum\GG J_k+2\ss
J=\AA\oo^2+2\AA\oo +\sum\AA_k\oo^2+\sum\ww J_k+2\ss J$. But
$\AA\oo^2+\sum\AA_k\oo^2+2\AA\oo=\GG (\la\AA,\AA_1,\ldots,\AA_s\ra )$ (we have
$\ss(\la\AA,\AA_1,\ldots,\AA_s\ra )=\nn (\la \AA,\AA_1,\ldots,\AA_s\ra
)=\AA\oo$). But $\ww (\la\AA,\AA_1,\ldots,\AA_s\ra )=\sum\AA\1\dd
(\AA\AA_k)+2\AA\oo$. (See [OM, p. 280]. We have $\AA\dd (\AA_k/\AA )=\AA\1\dd
(\AA\AA_k)$.) So $\GG J=\GG (\la\AA,\AA_1,\ldots,\AA_s\ra )+\sum\ww J_k+2\ss
J=\AA\oo^2+\sum\AA\1\dd (\AA\AA_k)+2\AA\oo+\sum\ww J_k+2\ss J=\AA\oo^2+\sum\ww
J_k+\sum\AA\1\dd (\AA\AA_k)+2\ss J$. (Recall, $2\AA\oo\sbq 2\ss J$.)
Let $\ww =\sum\ww J_k+\sum\AA\1\dd (\AA\AA_k)+2\ss J$. We have $\GG
J=\AA\oo^2+\ww$ and $2\ss J\sbq\ww$. By Lemma 2.10 in order to prove that $\ww
=\ww J$ we still need to prove that $\ww =2\ss J$ or $\od\AA +\od\ww$ is
odd. If $\ww\neq 2\ss J$, i.e. $\ww\sp 2\ss J$,
then $\ww =\AA\1\dd (\AA\AA_k)$ or $\ww =\ww J_k$ for some
$k$. Suppose that $\ww =\ww J_k$. We cannot have $\ww J_k=2\ss J_k\sbq
2\ss J$. So $\od\AA_k+\od\ww J_k$ is odd which implies that
$\od\AA+\od\ww J_k$ is odd unless $\od (\AA\AA_k)$ is odd. But this
would imply that $\AA_k\oo =\AA\1\dd (\AA\AA_k)\sbq\ww =\ww J_k$ so
$\ww J_k=\AA_k\oo$ which contradicts the fact that $\od\AA_k+\od\ww
J_k$ is odd. Finally if $\ww =\AA\1\dd (\AA\AA_k)$ then $\od\AA
+\od\ww=\od\dd (\AA\AA_k)$ is odd unless $\AA\AA_k\in\D\fs$. (If
$\a\in\ff$ has odd order then $\dd (\a )=\a\oo$ has odd order. If
$\od\a$ is even then $\od\dd (\a )=\od\a +d(\a )\ev d(\a )\m2$ is even
iff $d(\a )=2e$ i.e. iff $\a\in\D\fs$.) But this implies
that $\dd (\AA\AA_k)=4\AA\AA_k\oo$ i.e. $\ww =\AA\1\dd
(\AA\AA_k)=4\AA_k\oo\sb 2\ss J$. Contradiction. \qed 

\blm Suppose that $\nn L_k=\nn L^{\ss_k}$, $\nn L_{k+1}=\nn L^{\ss_{k+1}}$
and $\AA_k$ and $\AA_{k+1}$ are norm generators for $L_k$ and $L_{k+1}$,
respectively. If $u_k+u_{k+1}$ is even, then 
$$\FF_k=\ss_k^{-2}\dd (\AA_k\AA_{k+1})+\AA_k\ss_k^{-2}\ww
L_{(k+1)}^*+\AA_{k+1}\ww L_{(k)}^\* +2\p^{(u_k+u_{k+1})/2-r_k}.$$ 
\elm
\pf We have $L^{\ss_k}=\ss_kL_{(k)}^\*\pp L_{(k+1)}^*$ and
$L^{\ss_{k+1}}=\ss_{k+1}L_{(k)}^\*\pp L_{(k+1)}^*$. Now $L_{k+1}\sbq
L_{k+1}^*\sbq L^{\ss_{k+1}}$ and $L_k\sbq\ss_kL_{(k)}^\*\sbq L^{\ss_k}$. Thus
$\AA_{k+1}$ is norm generator for $L_{(k+1)}^*$ and for $L^{\ss_{k+1}}$ and
$\AA_k$ is a norm generator for $\ss_kL_{(k)}^\*$ and for $L^{\ss_k}$. Also
$\pi^{2(r_{k+1}-r_k)}\AA_k$ is a norm generator for
$\ss_{k+1}L_{(k)}^\*$. By Lemma 2.11 we get $\ww_k=\AA_k\1\dd
(\AA_k\AA_{k+1})+\ww (\ss_kL_{(k)}^\*)+\ww L_{(k+1)}^*+2\ss_k=\AA_k\1\dd
(\AA_k\AA_{k+1})+\ss_k^2\ww L_{(k)}^\*+\ww L_{(k+1)}^*+2\ss_k$ and
$\ww_{k+1}=\AA_{k+1}\1\dd (\pi^{2(r_{k+1}-r_k)}\AA_k\AA_{k+1})+\ww
(\ss_{k+1}L_{(k)}^\*)+\ww
L_{(k+1)}^*+2\ss_{k+1}=\AA_{k+1}\1\ss_{k+1}^2\ss_k^{-2}\dd
(\AA_k\AA_{k+1})+\ss_{k+1}^2\ww L_{(k)}^\*+\ww
L_{(k+1)}^*+2\ss_{k+1}$. 

By [OM, 93:26] we have $\ss_k^2\FF_k=\dd
(\AA_k\AA_{k+1})+\AA_{k+1}\ww_k+
\AA_k\ww_{k+1}+2\p^{(u_k+u_{k+1})/2+r_k}=\dd
(\AA_k\AA_{k+1})+\AA_k\1\AA_{k+1}\dd
(\AA_k\AA_{k+1})+\AA_{k+1}\ss_k^2\ww L_{(k)}^\*+\AA_{k+1}\ww
L_{(k+1)}^*+2\AA_{k+1}\ss_k+\\
\AA_k\AA_{k+1}\1\ss_{k+1}^2\ss_k^{-2}\dd
(\AA_k\AA_{k+1})+\AA_k\ss_{k+1}^2\ww L_{(k)}^\*+\AA_k\ww
L_{(k+1)}^*+2\AA_k\ss_{k+1}+2\p^{(u_k+u_{k+1})/2+r_k}$. 

But $\AA_k\oo\spq\AA_{k+1}\oo$ and $\AA_k\ss_k^{-2}\sbq
\AA_{k+1}\ss_{k+1}^{-2}$ ([OM, 93:25]) so $u_k\leq u_{k+1}$ and
$u_k-2r_k\geq u_{k+1}-2r_{k+1}$. Thus $\AA_k\1\AA_{k+1}\dd
(\AA_k\AA_{k+1}),\AA_k\AA_{k+1}\1\ss_{k+1}^2\ss_k^{-2}\dd
(\AA_k\AA_{k+1})\sbq \dd (\AA_k\AA_{k+1})$. Also $\AA_k\ss_{k+1}^2\ww
L_{(k)}^\*\sbq\AA_{k+1}\ss_k^2\ww L_{(k)}^\*$ and $\AA_{k+1}\ww
L_{(k+1)}^*\sbq\AA_k\ww L_{(k+1)}^*$. Also
$\od\AA_{k+1}\ss_k=u_{k+1}+r_k\geq (u_k+u_{k+1})/2+r_k$ (we have
$u_{k+1}\geq u_k$) and $\od\AA_k\ss_{k+1}=u_k+r_{k+1}\geq
(u_k+u_{k+1})/2+r_k$ (we have $u_k-2r_k\geq u_{k+1}-2r_{k+1}$). Hence
$2\AA_{k+1}\ss_k,2\AA_k\ss_{k+1}\sbq 2\p^{(u_k+u_{k+1})/2+r_k}$. 

By removing all unnecessary terms (which are included in others) we
get $\ss_k^2\FF_k=\dd (\AA_k\AA_{k+1})+\AA_{k+1}\ss_k^2\ww L_{(k)}^\*+
\AA_k\ww L_{(k+1)}^*+2\p^{(u_k+u_{k+1})/2+r_k}$. When we divide by
$\ss_k^2$ we get the desired result. \qed 

Suppose $L\ap\[ a_1,\ldots,a_n\]$ relative to the good BONG
$x_1,\ldots,x_n$. Let $L=L^1\pp\ldots\pp L^m$ be a maximal norm splitting
with all the binary components improper s.t. $x_1,\ldots,x_n$ is obtained
by putting together the BONGs of $L^1,\ldots,L^m$. We choose the Jordan
decomposition $L=L_1\pp\ldots\pp L_t$ with components obtained by putting
together the $L^j$'s of the same scale (see also the proof of [B, Lemma
4.7]). So the $L^j$'s with $\ss L^j=\ss_k$ make a maximal norm
splitting for $L_k$, those with $\ss L^j\sbq\ss_k$ a maximal norm
splitting for $L_{(k)}$ and those with $\ss L^j\sp\ss_k$ a maximal
norm splitting for $L^*_{(k+1)}$. By putting together the BONGs of the
components of these maximal norm splittings we get good BONGs for
$L_k$, $L_{(k)}$ and $L^*_{(k+1)}$. It follows that $L_k=\[
x_{n_{k-1}+1},\ldots,x_{n_k}\]$, $L_{(k)}=\[ x_1,\ldots,x_{n_k}\]$ and
$L^*_{(k+1)}=\[ x_{n_k+1},\ldots,x_n\]$. Also $\nn L_k=\nn
L^{\ss_k}$. (For any $L^j$ with $\ss L^j=\ss_k$ we have $L^j\sbq
L_k\sbq L^{\ss_k}$ and $\nn L^j=\nn L^{\ss L^j}=\nn L^{\ss_k}$.)

\blm (i) For any $n_{k-1}+1\leq i\leq n_k$ we have $R_i=u_k$ if $i\ev
n_{k-1}+1\m2$ and $R_i=2r_k-u_k$ if $i\ev n_{k-1}\m2$. 

(ii) For any $n_{k-1}+1\leq i\leq n_k$ we have $R_i=u_k$ if $i\ev
n_k+1\m2$ and $R_i=2r_k-u_k$ if $i\ev n_k\m2$. 

(iii) $\pm a_{n_{k-1}+1}$ and $\pm\pi^{2u_k-2r_k}a_{n_k}$ are norm
generators for $L_k$ and for $L^{\ss_k}$.  
\elm
\pf If $L_k$ is improper then $\dim L_k$ is even so $n_{k-1}\ev
n_k\m2$. Also the sequence $R_{n_{k-1}+1},\ldots,R_{n_k}$ is
$u_k,2r_k-u_k,\ldots,u_k,2r_k-u_k$ so we get both (i) and (ii). If $L_k$
is proper then $u_k=r_k$ and the sequence $R_{n_{k-1}+1},\ldots,R_{n_k}$
is $r_k,\ldots,r_k$. But $u_k=r_k$ so $r_k=u_k=2r_k-u_k$ and again we get
both (i) and (ii). 

(iii) We have $L_k\ap\[ a_{n_{k-1}+1},\ldots,a_{n_k}\]$ so
$a_{n_{k-1}+1}$ is a norm generator for $L_k$. We have
$L^\*_k\ap\[ a_{n_k}\1,\ldots,a_{n_{k-1}}\1\]$ so $a_{n_k}\1$ is a norm
generator for $L^\*_k=\p^{-r_k}L_k$. Therefore $\pi^{2r_k}a_{n_k}\1$
is a norm generator for $L_k$. But $\od a_{n_k}=2r_k-u_k$ so
$\pi^{2u_k-4r_k}a_{n_k}$ differs from $a_{n_k}\1$ by the square of a
unit. Since $\pi^{2r_k}a_{n_k}\1$ is a norm generator for $L_k$ so
is $\pi^{2r_k}\pi^{2u_k-4r_k}a_{n_k}=\pi^{2u_k-2r_k}a_{n_k}$. Since
$\GG L_k$ is an additive group $-a_{n_{k-1}+1}$ and
$-\pi^{2u_k-2r_k}a_{n_k}$ will also be norm genrators for $L_k$. We have
$L_k\sbq L^{\ss_k}$ and $\nn L_k=\nn L^{\ss_k}$ so $\pm
a_{n_{k-1}+1}$ and $\pm\pi^{2u_k-2r_k}a_{n_k}$ are norm generators for $L^{\ss
L_k}$ as well. \qed

We want now to find relations between $\a_i$'s and the O'Meara's invariants
$\ww_k$ and $\FF_k$. In particular, this will prove that $\a_i$'s are
invariants of the lattice $L$ i.e. they do not depend on the choice of the
BONG of $L$. 

\blm $\od\ww L=\min\{ R_1+\a_1,R_1+e\}$. (If $n=1$ we ignore
$R_1+\a_1$.) 

If moreover $L_1$ is not unary then $\od\ww L=R_1+\a_1$. 
\elm
\pf Note that if $L_1$ is not unary, in particular if $L^1$ is binary,
then $R_1=u_1\geq 2r_1-u_1=R_2$ so $\a_1\leq (R_2-R_1)/2+e\leq
e$. Hence $\min\{ R_1+\a_1,R_1+e\} =R_1+\a_1$ and so the two statements
of the lemma are equivalent. 

We use induction on $m$, the number of components in the maximal norm
splitting we fixed for $L$. Suppose first that $m=1$. If $L=L^1$ is
unary then $\ww L=2\ss L=2\p^{R_1}$ so $\od\ww L^1=R_1+e$, as
claimed. If $L=L^1$ is binary and so improper modular then we may
assume that it is unimodular since the statement is invariant upon
scaling. Hence $R_1+R_2=0$ and $R_1=\od\nn L>\od\ss L=0$. Now $a_1\in
Q(L)$ is a norm generator. Thus by [OM, 93:10] there is $b\in\ww L$
s.t. $L\ap A(a_1,b)$. Also if $\ww L\sp 2\ss L=2\oo$ then $\ww
L=b\oo$. Suppose first that $\ww L=2\oo$. Then $b\in 2\oo$ so $\od
b\geq e$. Thus $d(-a_1a_2)=d(-\det L)=d(1-a_1b)\geq\od a_1b\geq
R_1+e$ so $R_2-R_1+d(-a_1a_2)=-2R_1+d(-a_1a_2)\geq -R_1+e$. On the
other hand $(R_2-R_1)/2+e=-R_1+e$ so $\a_1=\min\{
(R_2-R_1)/2+e,R_2-R_1+d(-a_1a_2)\} =-R_1+e$. Thus
$\od\ww L=e=R_1+\a_1$. If $\ww L\sp 2\ss L=2\oo$ then $\ww L=b\oo$
and $\od a_1+\od b$ is odd. Also $\od a_1=\od\nn L\leq\od 2\ss L=e$
and $\od b=\od\ww L<\od 2\ss L=e$. It follows that $\od a_1b<2e$ and
it is odd. Hence $d(-a_1a_2)=d(1-a_1b)=\od a_1b=R_1+\od b$ so
$R_2-R_1+d(-a_1a_2)=-2R_1+d(-a_1a_2)=-R_1+\od b$. Also
$(R_2-R_1)/2+e=-R_1+e>-R_1+\od b$. It follows that
$\a_1=-R_1+\od b=-R_1+\od\ww L$. So $\od\ww L=R_1+\a_1$. 

We now prove the induction step. We have $L=L^1\pp L'$, where
$L'=L^2\pp\ldots\pp L^m$. Let now $\AA$ and $\AA'$ be a norm generators for
$L^1$ and $L'$. We have $\nn L^1=\nn L$ so $\AA$ is also a norm
generator for $L$. By Lemma 2.11 we have $\ww L=\ww L^1+\ww
L'+\AA\1\dd (\AA\AA')$. ($\AA\1\dd (\AA\AA )=0$ and $2\ss L=2\ss L^1\sbq\ww
L^1$ can be ignored.) Since $\od\AA\1\dd (\AA\AA')=\od\AA'+d(\AA\AA')$ it
follows that $\od\ww L=\min\{\od\ww L^1,\od\ww L',\od\AA'+d(\AA\AA')\}$. 

If $L^1$ is unary then $R_1\leq R_2$, $L^1\ap\[ a_1\]$ and $L'\ap\[
a_2,\ldots,a_n\]$. We take $\AA =a_1$ and $\AA'=-a_2$. We have $\od\AA'=R_2$,
$\ww L^1=R_1+e$ and $\od\ww L'=\min\{ R_2+\a_1(L'),R_2+e\}$. It
follows that $\od\ww L=\min\{
R_1+e,R_2+\a_1(L'),R_2+e,R_2+d(-a_1a_2)\}$. Since $R_2+e\geq R_1+e$,
it can be removed. By Corollary 2.5 (ii) we have $\a_1=\min\{
(R_2-R_1)/2+e,R_2-R_1+d(-a_1a_2),R_2-R_1+\a_1(L')\}$. It follows that
$\min\{ R_1+\a_1,R_1+e\} =\min\{ (R_1+R_2)/2+e, R_2+d(-a_1a_2),
R_2+\a_1(L'),R_1+e\}$. But $R_2\geq R_1$ so $(R_1+R_2)/2+e\geq
R_1+e$. Thus $\min\{ R_1+\a_1,R_1+e\} =\min\{ R_2+d(-a_1a_2),
R_2+\a_1(L'), R_1+e\} =\od\ww L$. 

If $L^1$ is binary then $R_1\geq R_2$, $L^1\ap\[ a_1,a_2\]$ and
$L'\ap\[ a_3,\ldots,a_n\]$. We prove that $\od\ww L=R_1+\a_1$. We take
$\AA =\pi^{2u_1-2r_1}a_2$ and $\AA'=-a_3$. (See Lemma 2.13(iii).) We have
$\ww L^1=R_1+\a_1(L^1)$, $\ww L'=\min\{ R_3+\a_1(L'),R_3+e\}$ and
$\od\AA'+d(\AA\AA')=R_3+d(-a_2a_3)$. Thus $\od\ww L=\min\{ R_1+\a_1(L^1),
R_3+\a_1(L'),R_3+d(-a_2a_3),R_3+e\}$. But $e\geq
(R_2-R_1)/2+e\geq\a_1(L')$ so $R_3+e\geq R_3+\a_1(L^1)$ and so $R_3+e$
can be removed. On the other hand $\a_1=\min\{ \a_1(L^1),
R_3-R_1+d(-a_2a_3), R_3-R_1+\a_1(L')\}$. (We have $\a_1(L^1)=\min\{
(R_2-R_1)/2+e,R_2-R_1+d(-a_1a_2)\}$ and, by Lemma 2.4(ii),
$R_3-R_1+\a_1(L')=R_3-R_1+\a_1(\[ a_3,\ldots,a_n\] )$ can replace all
$R_{j+1}-R_1+d(-a_ja_{j+1})$ with $j\geq 3$.) So $R_1+\a_1=\min\{
R_1+\a_1(L^1),R_3+d(-a_2a_3),R_3+\a_1(L')\} =\od\ww L$. \qed

\blm If $L_k$ is unary then $\ww_k=\ss_k(\FF_{k-1}+\FF_k+2\oo )$. (The
term $\FF_{k-1}$ is ignored if $k=1$ and $\FF_k$ is ignored if $k=t$.) 
\elm

\pf Since $L_k$ is unary we have $\ss_k=\AA_k\oo$ and $u_k=r_k$. Also
$\ww L_k=2\ss_k$. 

We have $L^{\ss_k}=(\perp_{j<k}\ss_k\ss_j\1 L_j)\perp L_k\pp
(\perp_{j>k}L_j)$. The first orthogonal sum is included in $\ss_{k-1}\1\ss_k
L^{\ss_{k-1}}$, while the last one is included in $L^{\ss_{k+1}}$. Hence
$L^{\ss_k}\sbq L_k+\ss_{k-1}\1\ss_kL^{\ss_{k-1}}+L^{\ss_{k+1}}$. The reverse
inclusion follows from [OM, 93:24] so
$L^{\ss_k}=L_k+\ss_{k-1}\1\ss_kL^{\ss_{k-1}}+L^{\ss_{k+1}}$. Now $\AA_k$ is a
norm generator for both $L^{\ss_k}$ and $L_k$, $\pi^{2(r_k-r_{k-1})}\AA_{k-1}$
for $\ss_{k-1}\1\ss_kL^{\ss_{k-1}}$ and $\AA_{k+1}$ for $L^{\ss_{k+1}}$. By
Lemma 2.11 we have $\ww_k=\AA_k\1\dd
(\pi^{2(r_k-r_{k-1})}\AA_{k-1}\AA_k)+\AA_k\1\dd (\AA_k\AA_{k+1})+\ww L_k+\ww
(\ss_{k-1}\1\ss_kL^{\ss_{k-1}})+\ww
L^{\ss_{k+1}}+2\ss_k=\ss_{k-1}^{-2}\ss_k\dd (\AA_{k-1}\AA_k)+\ss_k\1\dd 
(\AA_k\AA_{k+1})+\ss_{k-1}^{-2}\ss_k^2\ww_{k-1}+\ww_{k+1}+2\ss_k$. (We
ignore $\AA_k\dd (\AA_k\AA_k)=0$.)

If  $u_k+u_{k+1}$ is even then by [OM, 93:26] we have $\ss_k^2\FF_k=\dd
(\AA_k\AA_{k+1})+\AA_{k+1}\ww_k+\AA_k\ww_{k+1}+2\p^{(u_k+u_{k+1})/2+r_k})$. This
formula also holds in the case when $u_k+u_{k+1}$ is odd if we drop
the last term. Indeed, in this case $\ss_k^2\FF_k=\AA_k\AA_{k+1}\oo$ but
$\od\AA_k\AA_{k+1}$ is odd so $\dd (\AA_k\AA_{k+1})=\AA_k\AA_{k+1}\oo$ and we
also have $\AA_{k+1}\ww_k,\AA_k\ww_{k+1}\sbq\AA_k\AA_{k+1}\oo$. It follows
that $\ss_k(\FF_{k-1}+\FF_k+2\oo )=\ss_k(\ss_{k-1}^{-2}\dd
(\AA_{k-1}\AA_k)+\ss_{k-1}^{-2}\AA_k\ww_{k-1}+
\ss_{k-1}^{-2}\AA_{k-1}\ww_k+2\p^{(u_{k-1}+u_k)/2-r_{k-1}}
+\ss_k^{-2}\dd (\AA_k\AA_{k+1})+\ss_k^{-2}\AA_{k+1}\ww_k+
\ss_k^{-2}\AA_k\ww_{k+1}+2\p^{(u_k+u_{k+1})/2-r_k}+2\oo )$. (If
$u_k+u_{k+1}$ is odd we ignore $2\p^{(u_k+u_{k+1})/2-r_k}$. If
$u_{k-1}+u_k$ is odd we ignore $2\p^{(u_{k-1}+u_k)/2-r_{k-1}}$.) But
$r_k=u_k$ so
$(u_{k-1}+u_k)/2-r_{k-1}=(u_{k-1}-2r_{k-1}+2r_k-u_k)/2\geq 0$ and
$(u_k+u_{k+1})/2-r_k=(u_{k+1}-u_k)/2\geq 0$. Hence
$2\p^{(u_{k-1}+u_k)/2-r_{k-1}},2\p^{(u_k+u_{k+1})/2-r_k}\sbq 2\oo$ so
these terms can be ignored. Thus $\ss_k(\FF_{k-1}+\FF_k+2\oo
)=\ss_{k-1}^{-2}\ss_k\dd
(\AA_{k-1}\AA_k)+\ss_{k-1}^{-2}\ss_k^2\ww_{k-1}+
\ss_{k-1}^{-2}\ss_k\AA_{k-1}\ww_k+\ss_k\1\dd
(\AA_k\AA_{k+1})+\ss_k\1\AA_{k+1}\ww_k+\ww_{k+1}+2\ss_k=
\ww_k+\ss_{k-1}^{-2}\ss_k\AA_{k-1}\ww_k+\ss_k\1\AA_{k+1}\ww_k=\ww_k$.
(We have $\ss_{k-1}^{-2}\ss_k\AA_{k-1}=
\ss_{k-1}^{-2}\AA_{k-1}(\ss_k^{-2}\AA_k)\1\sbq\oo$ and
$\ss_k\1\AA_{k+1}=\AA_k\1\AA_{k+1}\oo\sbq\oo$. ) \qed

\blm Let $1\leq i\leq n-1$. Then: 

(i) If $n_{k-1}<i<n_k$ for some $1\leq k\leq t$, then
$R_i+\a_i=\od\ww_k$ and $-R_{i+1}+\a_i=\od\ww^\*_{t-k}$. 

(ii) Suppose that $i=n_k$ for some $1\leq k\leq t-1$. If $R_{i+1}-R_i$ is even
or $\leq 2e$ then $\a_i=\od\FF_k$; otherwise $\a_i=(R_{i+1}-R_i)/2+e$,
$\od\FF_k=R_{i+1}-R_i=2\a_i-2e$ and both $\a_i$ and $\od\FF_k$ are $>2e$. 
\elm
\pf (i) Note that $R_i+R_{i+1}=u_k+2r_k-u_k=2r_k$. Thus if
$R_i+\a_i=\od\ww_k$ then
$-R_{i+1}+\a_i=\od\ww_k-2r_k=\od\ss_k^{-2}\ww_k=\od\ww^\*_{t-k}$ so it
is enough to prove the first part of the statement. Also
$R_{n_{k-1}+1}+R_{n_{k-1}+2}=R_{n_k-1}+R_{n_k}=2r_k$ and so
$R_{n_{k-1}+1}+\a_{n_{k-1}+1}=\ldots=R_{n_k-1}+\a_{n_k-1}$ by Corollary
2.3(i). Thus it is enough to prove our statement for only one value of
$n_{k-1}<i<n_k$, say $i=n_{k-1}+1$. 

We use induction on $t$. Note that if $k=1$ then $\od\ww_1=\od\ww
L=R_1+\a_1$ by Lemma 2.14 so we are done. In particular, (i) is true
when $t=1$. Suppose now that $t\geq 2$. We may assume that $k\geq
2$. We have $L^{\ss_k}=\ss_k L_{(k-1)}^\*\pp L^*_{(k)}$. Since
$i=n_{k-1}+1<n_k$ we have $R_{i-1}=2r_{k-1}-u_{k-1}$, $R_i=u_k$,
$R_{i+1}=2r_k-u_k$, $L_{(k-1)}\ap\[ a_1,\ldots,a_{i-1}\]$ and
$L^*_{(k)}\ap\[ a_i,\ldots,a_n\]$. Note that $R_i\geq R_{i+1}$. If $\AA$ and
$\BB$ are norm generators for  $L^*_{(k)}$ and $\ss_kL_{(k-1)}^\*$,
respectively, then $\nn L^*_{(k)}=\p^{R_i}=\p^{u_k}=\nn
L^{\ss_k}$. Therefore $\AA$ is also a norm generator for $L^{\ss_k}$. By
Lemma 2.11 we have $\ww_k=\ss L^{\ss_k}=\ww (\ss_kL_{(k-1)}^\* )+\ww
L^*_{(k)}+\AA\1\dd (\AA\BB )+2\ss_k$, which implies that
$\od\ww_k=\min\{\od\ww (\ss_kL_{(k-1)}^\* ),\od\ww L^*_{(k)},\od\BB +d(\AA\BB
),r_k+e\}$. Now $L_k$ is not unary so $\od\ww L^*_{(k)}=R_i+\a_1(L^*_{(k)})$
by Lemma 2.14. Also $L_{(k-1)}^\*\ap\[ a_{i-1}\1,\ldots,a_1\1\]$ and $\od
a_{i-1}\1 =-R_{i-1}$ so by Lemma 2.14 we have $\od\ww L_{(k-1)}^\* =\min\{
-R_{i-1}+\a_1(L_{(k-1)}^\* ), -R_{i-1}+e\} =\min\{
-R_{i-1}+\a_{i-2}(L_{(k-1)}), -R_{i-1}+e\}$. (We have $\a_1(L_{(k-1)}^\*
)=\a_{i-2}(L_{(k-1)})$ by 2.6.) It follows that $\od\ww (\ss_kL_{(k-1)}^\*
)=\min\{ 2r_k-R_{i-1}+\a_{i-2}(L_{(k-1)}),2r_k-R_{i-1}+e\}=\min\{
R_i+R_{i+1}-R_{i-1}+\a_{i-2}(L_{(k-1)}),R_i+R_{i+1}-R_{i-1}+e\}$. Now
$a_{i-1}\1$ is a norm generator for $L_{(k-1)}^\*$, so $\BB
:=\pi^{2r_k}a_{i-1}\1$ is a norm generator for $\ss_kL_{(k-1)}^\*$, and
$\AA :=-a_i$ is a norm generator for $L_{(k)}$. We get $\od\BB +d(\AA\BB
)=2r_k-R_{i-1}+d(-\pi^{2r_k}a_{i-1}\1a_i)=R_i+R_{i+1}-R_{i-1}+d(-a_{i-1}a_i)$.
Also  $r_k+e=(R_i+R_{i+1})/2+e$. Thus $\od\ww_k=\min\{ R_i+\a_1(L^*_{(k)}), 
R_i+R_{i+1}-R_{i-1}+d(-a_{j-1}a_j),
R_i+R_{i+1}-R_{i-1}+\a_{i-2}(L_{(k-1)}),R_i+R_{i+1}-R_{i-1}+e,
(R_i+R_{i+1})/2+e\}$. But $R_{i-1}\leq R_{i+1}\leq R_i$ so
$R_i+R_{i+1}-R_{i-1}+e\geq (R_i+R_{i+1})/2+e=R_i+(R_{i+1}-R_i)/2+e\geq
R_i+\a_1(L^*_{(k)})$ so the last two terms can be removed. By Lemma 2.4(i)
$R_{i+1}-R_{i-1}+\a_{i-2}(L_{(k-1)})=R_{i+1}-R_{i-1}+\a_{i-2}(\[
a_1,\ldots,a_{i-1}\] )$ replaces all the terms in the definition of
$\a_i$ with $1\leq j\leq i-2$, while by Lemma 2.1 $\a_1(L^*_{(k)})=\a_1(\[
a_i,\ldots,a_n\] )$ replaces $(R_{i+1}-R_i)/2+e$ and the terms with
$i\leq j<n$. Hence $\a_i=\min\{ R_{i+1}-R_{i-1}+\a_{i-2}(L_{(k-1)}),
R_{i+1}-R_{i-1}+d(-a_{i-1}a_i),\a_1(L^*_{(k)})\}$. It follows that
$R_i+\a_i=\min\{ R_i+R_{i+1}-R_{i-1}+\a_{i-2}(L_{(k-1)}),
R_i+R_{i+1}-R_{i-1}+d(-a_{i-1}a_i),R_i+\a_1(L^*_{(k)})\} =\od\ww_k$. 

(ii) Since $i=n_k$ we have $R_i=2r_k-u_k$, $R_{i+1}=u_{k+1}$,
$L_{(k)}\ap\[ a_1,\ldots,a_i\]$ and $L^*_{(k+1)}\ap\[
a_{i+1},\ldots,a_n\]$. We have $R_{i+1}-R_i=u_k+u_{k+1}-2r_k$ so
$R_{i+1}-R_i$ is even iff $u_k+u_{k+1}$ is even. If $u_k+u_{k+1}$ and
$R_{i+1}-R_i$ are odd then $\FF_k=\ss_k^{-2}\AA_k\AA_{k+1}$ so
$\od\FF_k=u_k+u_{k+1}-2r_k=R_{i+1}-R_i$. If $R_{i+1}-R_i<2e$ then
$\a_i=R_{i+1}-R_i=\od\FF_k$, while if $R_{i+1}-R_i>2e$ then
$\a_i=(R_{i+1}-R_i)/2+e$. (See Lemmas 2.7(ii) and (iii).)

Suppose now that $R_{i+1}-R_i$ is even. By Lemma 2.12 we have
$\FF_k=\ss_k^{-2}\dd (\AA_k\AA_{k+1})+\AA_k\ss_k^{-2}\ww
L_{(k+1)}^*+\AA_{k+1}\ww L_{(k)}^\* +2\p^{(u_k+u_{k+1})/2-r_k}$. We
take $\AA_k=\pi^{2u_k-2r_k}a_i$ and $\AA_{k+1}=-a_{i+1}$. (See Lemma
2.13(iii).) Thus $\od\dd (\AA_k\AA_{k+1})=\od
(\AA_k\AA_{k+1})+d(\AA_k\AA_{k+1})=u_k+u_{k+1}+d(-a_ia_{i+1})$ and so
$\od\ss_k^{-2}\dd (\AA_k\AA_{k+1})=-2r_k+u_k+u_{k+1}+d(-a_ia_{i+1})=
R_{i+1}-R_i+d(-a_ia_{i+1})$. By Lemma 2.14 we have $\od\ww
L_{(k+1)}^*=\min\{ R_{i+1}+\a_1(L_{k+1}^*), R_{i+1}+e\}$. Since
$\od\AA_k\ss_k^{-2}=u_k-2r_k=-R_i$, we get $\od (\AA_k\ss_k^{-2}\ww
L_{(k+1)}^*)=\min\{ R_{i+1}-R_i+\a_1(L_{(k+1)}^*),R_{i+1}-R_i+e\}$. We
have $L_{(k)}^\*\ap\[ a_i\1,\ldots,a_1\1\]$, so $\od\ww L_{(k)}^\*
=\min\{ -R_i+\a_1(L_{(k)}^\* ),-R_i+e\}$. Since
$\a_1(L_{(k)}^\* )=\a_{i-1}(L_{(k)})$ (see 2.6) and
$\od\AA_{k+1}=u_{k+1}=R_{i+1}$ we have $\od (\AA_{k+1}\ww L_{(k)}^\*
)=\min\{ R_{i+1}-R_i+\a_{i-1}(L_{(k)}), R_{i+1}-R_i+e\}$. Finally $\od
2\p^{(u_k+u_{k+1})/2-r_k}=(u_k+u_{k+1})/2-r_k+e=(R_{i+1}-R_i)/2+e$.
Thus $\od\FF_k=\min\{ (R_{i+1}-R_i)/2+e, R_{i+1}-R_i+d(-a_ia_{i+1}),
R_{i+1}-R_i+\a_{i-1}(L_{(k)}), R_{i+1}-R_i+\a_1(L_{(k+1)}^*),
R_{i+1}-R_i+e\}$. But $R_i=2r_k-u_k\leq u_k\leq u_{k+1}=R_{i+1}$ so
$R_{i+1}-R_i+e\geq (R_{i+1}-R_i)/2+e$ so it can be ignored. So
$\od\FF_k=\min\{ (R_{i+1}-R_i)/2+e, R_{i+1}-R_i+d(-a_ia_{i+1}),
R_{i+1}-R_i+\a_{i-1}(L_{(k)}), R_{i+1}-R_i+\a_1(L_{(k+1)}^*)\}$,
which, by Corollary 2.5(ii), is equal to $\a_i$. (Recall,
$L_{(k)}\ap\[ a_1,\ldots,a_i\]$ and $L_{(k+1)}^*\ap\[
a_{i+1},\ldots,a_n\]$.)\qed

\bco (i) If $L_k$ is not unary and $i=n_{k-1}+1$ or $n_k-1$ then
$\od\AA_k\1\ww_k=\a_i$. 

(ii) If $L_k$ is unary and $i=n_k$ then
$\od\AA_k\1\ww_k=\min\{\a_{i-1},\a_i,e\}$. (We ignore $\a_{i-1}$ if
$i=1$, and $\a_i$ if $i=n$.) 
\eco
\pf (i) In both cases when $i=n_{k-1}+1$ or $n_k-1$ we have
$R_i=u_k=\od\AA_k$. Hence $\od\ww_k=R_i+\a_i=\od\AA_k+\a_i$ so
$\od\AA_k\1\ww_k=\a_i$. 

(ii) We have $\ss_k=\AA_k\oo$ and, by Lemma 2.15,
$\ww_k=\ss_k(\FF_{k-1}+\FF_k+2\oo )$ so
$\AA_k\1\ww_k=\FF_{k-1}+\FF_k+2\oo$. Thus
$\od\AA_k\1\ww_k=\min\{\od\FF_{k-1},\od\FF_k,e\}$ and we have to prove
that it is equal to $\min\{\a_{i-1},\a_i,e\}$. Now $i-1=n_k-1=n_{k-1}$
so, by Lemma 2.16(ii), we have either $\a_{i-1}=\od\FF_{k-1}$ or
$\a_{i-1},\od\FF_{k-1}>2e$. But if $\a_{i-1},\od\FF_{k-1}>2e>e$ then
they can be ignored in $\min\{\a_{i-1},\a_i,e\}$
and $\min\{\od\FF_{k-1},\od\FF_k,e\}$, respectively. Similarly either
$\a_i=\od\FF_k$ or $\a_i,\od\FF_k$ are both $>2e>e$ so they can be
ignored. Thus $\min\{\a_{i-1},\a_i,e\}
=\min\{\od\FF_{k-1},\od\FF_k,e\}$. \qed

\section{Main theorem}

In this section we state and prove the main result of this paper, the
classification of integral lattices over dyadic local fields in terms of good
BONGs. It is well known that this problem was first solved by O'Meara in [OM,
Theorem 93:28]. Since our proof uses O'Meara's result we first state Theorem
93:28. 

Throughout this section $L,K$ are two lattices with $L\ap\[ a_1,\ldots,a_n\]$
and $K\ap\[ b_1,\ldots,b_n\]$ relative to good BONGs. In terms of Jordan
decompositions we write $L=L_1\pp\ldots\pp L_t$ and $K=K_1\pp\ldots\pp
K_{t'}$. Let $\ss_k=\ss L_k$, $\ss'_k=\ss K_k$, $\GG_k=\GG L^{\ss_k}$,
$\GG'_k=\GG K^{\ss'_k}$, $\ww_k=\ww L^{\ss_k}$, $\ww'_k=\ww K^{\ss'_k}$,
$\FF_k=\FF_k(L)$ and $\FF'_k=\FF_k(K)$. Let $\AA_k$ and $\BB_k$ be norm
generators for $L^{\ss_k}$ and $K^{\ss'_k}$, respectively. We say that $L$ and
$K$ are of the same fundamental type if
$$t=t',~\dim L_k=\dim K_k,~\ss_k=\ss'_k,~\GG_k=\GG'_k$$
for $1\leq k\leq t$. These conditions are equivalent to
$$t=t',~\dim L_k=\dim
K_k,~\ss_k=\ss'_k,~\ww_k=\ww'_k,~\AA_k\ap\BB_k\pmod{\ww_k}$$
for $1\leq k\leq t$. We now state O'Meara's Theorem 93:28.
\\

\no{\bf Theorem 93:28} {\it Let $L,K$ be lattices with the same fundamental
type such that $FL\ap FK$. Let $L_{(1)}\sb\cdots\sb L_{(t)}$ and
$K_{(1)}\sb\cdots\sb K_{(t)}$ be Jordan chains for $L$ and $K$. Then $L\ap K$
if and only if the following conditions hold for $1\leq i\leq t-1$

(i) $\det L_{(k)}/\det K_{(k)}\ap 1\pmod{\FF_k}$.

(ii) $FL_{(k)}\rep FK_{(k)}\perp [\AA_{k+1}]$ when $\FF_k\sb
4\AA_{k+1}\ww_{k+1}\1$. 

(iii) $FL_{(k)}\rep FK_{(k)}\perp [\AA_k]$ when $\FF_k\sb 4\AA_k\ww_k\1$.} 
\\

We state now our main result.

\btm Let $L,K$ be two lattices with $FL\ap FK$ and let $L\ap\[
a_1,\ldots,a_n\]$ and $K\ap\[ b_1,\ldots,b_n\]$ relative to good BONGs. Let
$R_i=R_i(L)=\od a_i$, $S_i=R_i(K)=\od b_i$, $\a_i=\a_i(L)$ and
$\b_i=\a_i(K)$. Then $L\ap K$ iff: 

(i) $R_i=S_i$ for $1\leq i\leq n$

(ii) $\a_i=\b_i$ for $1\leq i\leq n-1$

(iii) $d(a_1\cdots a_ib_1\cdots b_i)\geq\a_i$ for $1\leq i\leq n-1$

(iv) $[b_1,\ldots,b_{i-1}]\rep [a_1,\ldots,a_i]$ for any $1<i<n$
s.t. $\a_{i-1}+\a_i>2e$. 
\etm
\pf Condition 3.1(i) is equivalent to $t=t'$, $\dim L_k=\dim K_k$,
$\ss_k=\ss'_k$ and $\nn L^{\ss_k}=\nn K^{\ss_k}$
i.e. $\AA_k\oo=\BB_k\oo$. (See [B, Lemma 4.7].) Suppose this happens. Denote
as before $n_k=\dim L_{(k)}=\dim K_{(k)}$, $\p^{r_k}=\ss_k$ and $\p^{u_k}=\nn
L^{\ss_k}=\AA_k\oo$. 

As in the previous section, we choose a Jordan splitting of $L$ such that
$L_k\ap\[ a_{n_{k-1}+1},\ldots,a_{n_k}\]$. Hence for any $1\leq k\leq n$,
$\AA_k$ can be either $\pm a_{n_{k-1}+1}$ or $\pm\pi^{2u_k-2r_k}a_{n_k}$. We
choose a Jordan splitting for $K$ with the same property. 

Assuming that 3.1(i) holds, Lemma 2.16 and Corollary 2.17(ii) imply that
3.1(ii) is equivalent to $\ww_k=\ww'_k$ for $1\leq k\leq t$ and
$\FF_k=\FF'_k$ for $1\leq k\leq t-1$. 

From here the proof of Theorem 3.1 consists of two steps:

1. Assuming that 3.1(i) and (ii) hold, we prove that condition 3.1(iii) is
equivalent to $\AA_k\ap\BB_k\pmod{\ww_k}$ for any $1\leq k\leq t$ and condition
93:28(i). 

2. Assuming that 3.1(i)-(iii) hold, we prove that condition 3.1(iv) is
equivalent to conditions 93:28(ii) and (iii). 
\\

\blm Suppose that $L,K$ satisfy conditions 3.1(i) and 3.1(ii). If
$R_{i-1}=R_{i+1}$ for some $1<i<n$ then: 

(i) If 3.1(iii) holds at $i-2$ or $i-2=0$ then 3.1(iii) holds at $i$. 

(ii) If 3.1(iii) holds at $i+1$ or $i+1=n$ then 3.1(iii) holds at $i-1$. 
\elm
\pf (i) We have $d(a_1\cdots a_ib_1\cdots b_i)\geq\min\{
d(a_1\cdots a_{i-2}b_1\cdots b_{i-2}),d(-a_{i-1}a_i),d(-b_{i-1}b_i)\}$. (If
$i-2=0$ we ignore $d(a_1\cdots a_{i-2}b_1\cdots b_{i-2})$.) But
$d(a_1\cdots a_{i-2}b_1\cdots b_{i-2})\geq\a_{i-2}\geq
R_{i-1}-R_{i+1}+\a_i=\a_i$. (We have $-R_{i-1}+\a_{i-2}\geq
-R_{i+1}+\a_i$.) Also
$d(-a_{i-1}a_i)=R_{i+1}-R_{i-1}+d(-a_{i-1}a_i)\geq\a_i$. Similarly
$d(-b_{i-1}b_i)\geq\a_i$. Hence
$d(a_1\cdots a_ib_1\cdots b_i)\geq\a_i$. 

(ii) is similar. This time $R_{i+1}+\a_{i+1}\geq R_{i-1}+\a_{i-1}$ so
$d(a_1\cdots a_{i+1},b_1\cdots b_{i+1})\geq\a_{i+1}\geq
R_{i-1}-R_{i+1}+\a_{i-1}=\a_{i-1}$. (If $i+1=n$ then
$d(a_1\cdots a_nb_1\cdots b_n)=\j >\a_{n-2}$) Also
$d(-a_ia_{i+1})=R_{i+1}-R_{i-1}+d(-a_ia_{i+1})\geq\a_{i-1}$ and similarly
$d(-b_ib_{i+1})\geq\a_{i-1}$. \qed 

\blm Assuming that 3.1(i) and (ii) hold, condition 3.1(iii) is equivalent to
$\AA_k\ap\BB_k\pmod{\ww_k}$ for any $1\leq k\leq t$ and condition 93:28(i). 
\elm
\pf We have $L_{(k)}\ap\[ a_1,\ldots,a_{n_k}\]$ and $K_{(k)}\ap\[
b_1,\ldots,b_{n_k}\]$. Hence $\det L_{(k)}=a_1\cdots a_{n_k}$ and
$\det K_{(k)}=b_1\cdots b_{n_k}$. Since the two determinants have the
same order, $R_1+\cdots +R_{n_k}$, the condition $\det L_{(k)}/\det
K_{(k)}\ap 1(\mo\FF_k)$ is equivalent to $d(a_1\cdots a_{n_k}b_1\cdots
b_{n_k})\geq\od\FF_k$. Let $i=n_k$. We claim that $d(a_1\cdots
a_ib_1\cdots b_i)\geq\od\FF_k$ is equivalent to $d(a_1\cdots
a_ib_1\cdots b_i)\geq\a_i$. By Lemma 2.16(ii) we have either
$\a_i=\od\FF_k$ or $\a_i,\od\FF_k>2e$. In the first case our claim is
obvious and in the second both $d(a_1\cdots a_ib_1\cdots
b_i)\geq\od\FF_k$ and $d(a_1\cdots a_ib_1\cdots b_i)\geq\a_i$ are
equivalent to $a_1\cdots a_ib_1\cdots b_i\in\fs$. 

Thus condition 3.1(iii) at indices $i=n_k$ with $1\leq k\leq t-1$ is
equivalent to 93:28(i). Assume these equivalent conditions hold. We want to
prove that condition $\AA_k\ap\BB_k\pmod{\ww_k}$ at indices $1\leq k\leq t$
s.t. $L_k$ is not unary is equivalent to condition 3.1(iii) at $i=n_{k-1}+1$,
while if $L_k$ is unary then it holds unconditionally. 

Note that $\AA_k\ap\BB_k(\mo\ww_k)$ is equivalent to $\BB_k/\AA_k\ap
1(\mo\AA_k\1\ww_k)$ i.e. to
$d(\AA_k\BB_k)=d(\BB_k/\AA_k)\geq\od\AA_k\1\ww_k$. We will take
$\AA_k=a_{n_{k-1}+1}=a_i$ and $\BB_k=b_{n_{k-1}+1}=b_i$. So our condition is
equivalent to $d(a_ib_i)\geq\od\AA_k\1\ww_k$, where $i=n_{k-1}+1$. 

If $L_k$ is unary then $\od\AA_k\1\ww_k=\min\{\a_{i-1},\a_i,e\}$ by
Corollary 2.17(ii), where $i=n_{k-1}+1=n_k$. Since $i-1=n_{k-1}$ and
$i=n_k$, condition 3.1(iii) is satisfied for both. Thus $d(a_1\cdots
a_{i-1}b_1\cdots b_{i-1})\geq\a_{i-1}\geq\od\AA_k\1\ww_k$ and $d(a_1\cdots
a_ib_1\cdots b_i)\geq\a_i\geq\od\AA_k\1\ww_k$ so
$d(a_ib_i)\geq\od\AA_k\1\ww_k$. (If $k=1$ so $i=n_0+1=1$ we ignore $\a_{i-1}$
and we have $d(a_1b_1)\geq\a_1\geq\od\AA_1\1\ww_1$. If $k=t$ so $i=n_t=n$ we
ignore $\a_i$ and, since $a_1\cdots a_n=\det FM=\det FN=b_1\cdots b_n$ in
$\ff/\fs$, we get $d(a_nb_n)=d(a_1\cdots a_{n-1}b_1\cdots
b_{n-1})\geq\a_{n-1}\geq\od\AA_t\1\ww_t$.) Thus condition
$\AA_k\ap\BB_k\pmod{\ww_k}$ is superfluous when $L_k$ is unary. 

Suppose now that $L_k$ is not unary and let $i=n_{k-1}+1$. By
Corollary 2.17(i) we have $\od\AA_k\1\ww_k=\a_i$. We will prove that
$d(a_ib_i)\geq\od\AA_k\1\ww_k=\a_i$ is equivalent to the condition
3.1(iii) at $i$ i.e. to $d(a_1\cdots a_ib_1\cdots b_i)\geq\a_i$. If
$k=1$ so $i=n_0+1=1$ this is obvious. If $k>1$ so $i>1$ note that
$-R_i+\a_{i-1}\geq -R_{i+1}+\a_i$ and $R_i=u_k\geq 2r_k-u_k=R_{i+1}$ so
$\a_{i-1}\geq\a_i$. We have $i-1=n_{k-1}$ so $d(a_1\cdots a_{i-1}b_1\cdots
b_{i-1})\geq\a_{i-1}\geq\a_i$ and so $d(a_1b_i)\geq\a_i$ is equivalent to
$d(a_1\cdots a_ib_1\cdots b_i)\geq\a_i$ by domination principle. 

To complete the proof we show that 3.1(iii) is true if it is true for $i=n_k$,
where $1\leq k\leq t-1$, and for $i=n_{k-1}+1$, where $1\leq k\leq t$ and
$L_k$ is not unary. To do this we use Lemma 3.2. 

Let $1\leq k\leq t$. For any $n_{k-1}+1<i<n_k$ we have
$R_{i-1}=R_{i+1}$ (they are both $u_k$ or $2r_k-u_k$) so by Lemma
3.2(i) if 3.1(iii) holds for $i-2$ or $i-2=0$ it will also hold for
$i$. Thus, since 3.1(iii) is true for $n_{k-1}$ (or $n_{k-1}=0$ if
$k=1$), it will also be true by induction for any $n_{k-1}+2\leq i<n_k$
with $i\ev n_{k-1}\m2$. Similarly since 3.1(iii) is true at $n_{k-1}+1$,
it will also be true by induction for any $n_{k-1}+1\leq i<n_k$ with
$i\ev n_{k-1}+1\m2$. Hence 3.1(iii) holds for any $n_{k-1}<i<n_k$. Since
3.1(iii) also holds for any $i=n_k$ with $1\leq k\leq t-1$ it will hold
for any $1\leq i\leq n-1$. \qed

\blm If $1<i<n$ and $R_{i-1}=R_{i+1}$ then $\a_{i-1}+\a_i\leq 2e$. 
\elm
\pf We have $\a_{i-1}+\a_i\leq
(R_i-R_{i-1})/2+e+(R_{i+1}-R_i)/2+e=(R_{i+1}-R_{i-1})/2+2e$ so if
$R_{i-1}=R_{i+1}$ then $\a_{i-1}+\a_i\leq 2e$. \qed

\blm Let $V,W$ be two quadratic spaces over $F$. We have: 

(i) If $\dim V-\dim W=1$ and $\hh$ is a hyperbolic plane then $W\rep
V$ iff $V\rep W\pp\hh$. 

(ii) If $\dim V=\dim W$ and $a\in\ff$ then $W\rep V\pp [a]$ iff $V\rep
W\pp [a\det V\det W]$. 

(iii) If $\dim V=\dim W$, $a,b\in\ff$ and $(ab,\det V\det W)_\p =1$
(in particular, if $d(ab)+d(\det V\det W)>2e$) then $W\rep V\pp [a]$
iff $W\rep V\pp [b]$. 
\elm
\pf This is a direct consequence of [OM, 63:21]. For (iii) we also use the
fact that if $xy=zt$ then $[x,y]\ap [z,t]$ iff $z\rep [x,y]$, which in turn is
equivalent to $(xz,yz)_\p =1$. \qed

\blm Suppose that $L,K$ satisfy the conditions 3.1(i)-(iii) (or, equivalently,
they have the same fundamental type and they satisfy the condition
93:28(i)). Then: 

(i) If $\FF_k\sb 4\AA_k\ww_k\1$ and both $\AA_k$ and $\BB_k$ are norm
generators for $L^{\ss_k}$, then $FL_{(k)}\rep FK_{(k)}\pp [\AA_k]$ is
equivalent to $FL_{(k)}\rep FK_{(k)}\pp [\BB_k]$, and also to
$FK_{(k)}\rep FL_{(k)}\pp [\BB_k]$. 

(ii) If $\FF_k\sb 4\AA_{k+1}\ww_{k+1}\1$ and both $\AA_{k+1}$ and $\BB_{k+1}$
are norm generators for $L^{\ss_{k+1}}$ then $FL_{(k)}\rep FK_{(k)}\pp
[\AA_{k+1}]$ is equivalent to $FL_{(k)}\rep FK_{(k)}\pp [\BB_{k+1}]$ and also
to $FK_{(k)}\rep FL_{(k)}\pp [\BB_{k+1}]$. 
\elm
\pf (i) $\FF_k\sb 4\AA_k\ww_k\1$ is equivalent to
$\od\FF_k+\od\AA_k\1\ww_k>2e$. We have $\AA_k\1\ww_k\spq
2\AA_k\1\ss_k\spq 2\oo$ so $\od\AA_k\1\ww_k\leq e<\od\FF_k$. Since
$\AA_k,\BB_k$ are both norm generators	for $L^{\ss_k}$ we have
$d(\AA_k\BB_k)\geq\od\AA_k\1\ww_k$. Since also $d(\det L_{(k)}\det
K_{(k)})\geq\od\FF_k>\od\AA_k\1\ww_k$ we also have $d(\AA_k\BB_k\det
L_{(k)}\det K_{(k)})\geq\od\AA_k\1\ww_k$. Since $d(\det L_{(k)}\det
K_{(k)})+d(\AA_k\BB_k)\geq\od\FF_k+\od\AA_k\1\ww_k>2e$ we get by Lemma
3.5(iii) that $FL_{(k)}\rep FK_{(k)}\pp [\AA_k]$ iff $FL_{(k)}\rep
FK_{(k)}\pp [\BB_k]$. Similarly, since $d(\det L_{(k)}\det
K_{(k)})+d(\AA_k\BB_k\det L_{(k)}\det
K_{(k)})\geq\od\FF_k+\od\AA_k\1\ww_k>2e$, we have $FL_{(k)}\rep
FK_{(k)}\pp [\AA_k]$ iff $FL_{(k)}\rep FK_{(k)}\pp [\BB_k\det
L_{(k)}\det K_{(k)}]$ which, by Lemma 3.5(ii), is equivalent to
$FK_{(k)}\rep FL_{(k)}\pp [\BB_k]$. 

(ii) Same proof from (i) but with $\AA_k,\BB_k,\ww_k$ replaced by
$\AA_{k+1},\BB_{k+1},\ww_{k+1}$. \qed 

\blm Suppose that $L,K$ satisfy the conditions 3.1(i) - (iii). If $1\leq k\leq
t-1$ then: 

(i) If $\FF_k\sb 4\AA_k\ww_k\1$ then $FL_{(k)}\rep FK_{(k)}\pp [\AA_k]$
iff $[b_1,\ldots,b_{i-1}]\rep [a_1,\ldots,a_i]$, with $i=n_k$.

(ii) If $\FF_k\sb 4\AA_{k+1}\ww_{k+1}\1$ then $FL_{(k)}\rep
FK_{(k)}\pp [\AA_{k+1}]$ iff $[b_1,\ldots,b_{i-1}]\rep [a_1,\ldots,a_i]$,
with $i=n_k+1$. 
\elm
\pf (i) We take $\BB_k=-\pi^{2u_k-2r_k}b_i$ as a norm generator for
$K^{\ss_k}$, so for $L^{\ss_k}$. (See Lemma 2.13(iii).) By Lemma 3.6(i)
$FL_{(k)}\rep FK_{(k)}\pp [\AA_k]$ iff $FL_{(k)}\rep FK_{(k)}\pp
[\BB_k]$ i.e. iff $[a_1,\ldots,a_i]\rep [b_1,\ldots,b_i]\pp [-b_i]\ap
[b_1,\ldots,b_{i-1}]\pp\hh$. By Lemma 3.5(i) this is equivalent to
$[b_1,\ldots,b_{i-1}]\rep [a_1,\ldots,a_i]$. 

(ii) We take $\BB_{k+1}=a_i$ as a norm generator for $L^{\ss_{k+1}}$. By Lemma
3.6(ii) $FL_{(k)}\rep FK_{(k)}\pp [\AA_{k+1}]$ iff $FK_{(k)}\rep FL_{(k)}\pp
[\BB_{k+1}]$ i.e. iff $[b_1,\ldots,b_{i-1}]\rep [a_1,\ldots,a_{i-1}]\pp
[a_i]\ap [a_1,\ldots,a_i]$. \qed


\blm (i) If $i=n_k>n_{k-1}+1$ then $\a_{i-1}+\a_i>2e$ iff
$\FF_k\sb\AA_k\ww_k\1$.

(ii) If $i=n_k+1<n_{k+1}$ then $\a_{i-1}+\a_i>2e$ iff
$\FF_k\sb\AA_{k+1}\ww_{k+1}\1$. 

(iii) If $i=n_k=n_{k-1}+1$ then $\a_{i-1}+\a_i>2e$ iff
$\FF_k\sb\AA_k\ww_k\1$ or $\FF_{k-1}\sb\AA_k\ww_k\1$. 

(In (iii) we ignore the condition $\FF_k\sb\AA_k\ww_k\1$ if $k=t$ and
we ignore $\FF_{k-1}\sb\AA_k\ww_k\1$ if $k=1$.)
\elm 
\pf (i) Condition $\FF_k\sb 4\AA_k\ww_k\1$ is equivalent to
$\od\AA_k\1\ww_k+\od\FF_k>2e$. By Corollary 2.17(i) we have
$\od\AA_k\1\ww_k=\a_{i-1}$. By Lemma 2.16(ii) we have either
$\a_i=\od\FF_k$ or $\a_i,\od\FF_k>2e$. In the first case
$\od\AA_k\1\ww_k+\od\FF_k=\a_{i-1}+\a_i$ and in the second both
$\od\AA_k\1\ww_k+\od\FF_k>2e$ and $\a_{i-1}+\a_i>2e$ hold. In both
cases $\od\AA_k\1\ww_k+\od\FF_k>2e$ iff $\a_{i-1}+\a_i>2e$.

(ii) We have $\FF_k\sb 4\AA_{k+1}\ww_{k+1}\1$ iff
$\od\AA_{k+1}\1\ww_{k+1}+\od\FF_k>2e$. By Corollary 2.17(i)
$\od\AA_{k+1}\1\ww_{k+1}=\a_i$ and by Lemma 2.16(ii) $\od\FF_k$ and
$\a_{i-1}$ are either equal or they are both $>2e$. Thus
$\od\AA_{k+1}\1\ww_{k+1}+\od\FF_k>2e$ iff $\a_{i-1}+\a_i>2e$.

(iii) $\FF_{k-1}\sb 4\AA_k\ww_k\1$ and  $\FF_k\sb 4\AA_k\ww_k\1$ are
equivalent to $\od\FF_{k-1}+\od\AA_k\1\ww_k>2e$ resp.
$\od\FF_k+\od\AA_k\1\ww_k>2e$. By Corollary 2.17(ii) we have
$\od\AA_k\1\ww_k=\min\{\a_{i-1},\a_i,e\}\geq 0$. By Lemma 2.16(ii) we
have that $\od\FF_{k-1}=\a_{i-1}$ or $\od\FF_{k-1},\a_{i-1}>2e$ and
$\od\FF_k=\a_i$ or $\od\FF_k,\a_i>2e$. Therefore $\FF_{k-1}\sb
4\AA_k\ww_k\1$ and $\FF_k\sb 4\AA_k\ww_k\1$ are equivalent to
$\a_{i-1}+\min\{\a_{i-1},\a_i,e\} >2e$
resp. $\a_i+\min\{\a_{i-1},\a_i,e\} >2e$. Obviously either of them
implies $\a_{i-1}+\a_i>2e$. Conversely, suppose that
$\a_{i-1}+\a_i>2e$. If both $\a_{i-1}$ and $\a_i$ are $>e$ then we
have both $\a_{i-1}+\min\{\a_{i-1},\a_i,e\} >2e$ and
$\a_i+\min\{\a_{i-1},\a_i,e\} >2e$. Otherwise we have
$\min\{\a_{i-1},\a_i,e\} =\min\{\a_{i-1},\a_i\}$ and so
$\max\{\a_{i-1},\a_i\} +\min\{\a_{i-1},\a_i,e\} =\max\{\a_{i-1},\a_i\}
+\min\{\a_{i-1},\a_i\} =\a_{i-1}+\a_i>2e$, which implies that either
$\a_{i-1}+\min\{\a_{i-1},\a_i,e\} >2e$ or
$\a_i+\min\{\a_{i-1},\a_i,e\} >2e$. \qed

\blm Assuming that 3.1(i)-(iii) hold, condition 3.1(iv) is equivalent to
93:28(ii) and (iii). 
\elm
\pf Take $1<i<n$. If $n_{k-1}+1<i<n_k$ for some $1\leq k\leq t$ then
$R_{i-1}=R_{i+1}$, by Lemma 2.13, so, by Lemma 3.4, $\a_{i-1}+\a_i\leq 2e$,
which makes 3.1(iv) vacuous at $i$. Therefore we can restrict ouselves to
$i=n_k$ or $n_k+1$ for some $1\leq k\leq t-1$. We have three cases: 

1. $i=n_k$ and $\dim L_k>1$ i.e. $i=n_k>n_{k-1}+1$. By Lemma 3.8(i)
$\FF_k\sb 4\AA_k\ww_k\1$ is equivalent to $\a_{i-1}+\a_i>2e$. On the
other hand if $\FF_k\sb 4\AA_k\ww_k\1$ then $FL_{(k)}\rep FK_{(k)}\pp
[\AA_k]$ is equivalent to $[b_1,\ldots,b_{i-1}]\rep [a_1,\ldots,a_i]$
by Lemma 3.7(i). Therefore 3.1(iv) at index $i$ is equivalent to
93:28(iii) at index $k$. 

2. $i=n_k+1$ and $\dim L_{k+1}>1$ i.e. $i=n_k+1<n_{k+1}$. By Lemma
3.8(ii) $\FF_k\sb 4\AA_{k+1}\ww_{k+1}\1$ is equivalent to
$\a_{i-1}+\a_i>2e$. On the other hand $FL_{(k)}\rep FK_{(k)}\pp
[\AA_{k+1}]$ is equivalent to $[b_1,\ldots,b_{i-1}]\rep
[a_1,\ldots,a_i]$ by Lemma 3.16(ii). We will prove that 3.1(iv) at
index $i$ is equivalent to 93:28(ii) at index $k$. 

3. $i=n_k=n_{k-1}+1$ for some $1\leq k\leq t$. In this case $L_k$ is unary. We
will prove that the condition 3.1(iv) at index $i$ is equivalent to 93:28(iii)
at index $k$ and 93:28(ii) at index $k-1$. First note that if $k=t$ then
3.1(iv) is vacuous at $i=n_t=n$. On the other hand 93:28(iii) is vacuous at
index $k=t$. Also if $\FF_{t-1}\sb 4\AA_t\ww_t\1$ then, by Lemma 3.7(ii),
$FL_{(t-1)}\rep FK_{(t-1)}\pp [\AA_t]$ is equivalent to
$[b_1,\ldots,b_{n-1}]\rep [a_1,\ldots,a_n]$ (we have $i=n_{t-1}+1=n_t=n$). But
this follows from $[a_1,\ldots,a_n]\ap [b_1,\ldots,b_n]$. Thus 93:28(ii) is
superfluous at index $k-1=t-1$. Next we note that if $k=1$ then 3.1(iv) is
vacuous at $i=n_0+1=1$. On the other hand 93:28(ii) is vacuous at index
$k-1=0$. Also if $\FF_1\sb 4\AA_1\1\ww_1$ then $FL_{(1)}\rep
FK_{(1)}\pp [\AA_1]$ is equivalent, by Lemma 3.7(i), to $0\rep [a_1]$
(we have $i=n_1=1$). Here $0$ is not the scalar zero, but the zero
lattice, of dimension $0$, so $0\rep [a_1]$ holds trivially. Thus 93:28(iii)
is superfluous at $k=1$. 

Suppose now that $1<k<t$. By Lemma 3.8(iii) we have $\a_{i-1}+\a_i>2e$
iff $\FF_{k-1}\sb 4\AA_k\ww_k\1$ or $\FF_k\sb 4\AA_k\ww_k\1$.  To
complete the proof we note that if $\FF_{k-1}\sb 4\AA_k\ww_k\1$ then
$FL_{(k-1)}\rep FK_{(k-1)}\pp [\AA_k]$ is equivalent to
$[b_1,\ldots,b_{i-1}]\rep [a_1,\ldots,a_i]$ by Lemma 3.7(ii) (we have
$i=n_{k-1}+1$) and if $\FF_k\sb 4\AA_k\ww_k\1$ then $FL_{(k)}\rep
FK_{(k)}\pp [\AA_k]$ is equivalent to $[b_1,\ldots,b_{i-1}]\rep
[a_1,\ldots,a_i]$ by Lemma 3.7(i) (we have $i=n_k$). \qed

 
\section{The 2-adic case}

In this section we will assume that $F$ is 2-adic i.e. that $e=1$. 

In [OM, \S 93G] O'Meara gives a solution to the classification problem in the
2-adic case which only involves the Jordan invariants $t$, $\dim L_k$, $\ss_k$
and $\nn_k:=\nn L_k$. The invariants $\GG_k$ and $\ww_k$ are no longer
necessary since they can be written as $\GG_k=\nn_k$ and $\ww_k=2\ss_k$. A
similar phenomenon occurs when we use good BONGs instead of Jordan
decompositions. This time the invarians $\a_i$ are no longer necessary. 

\blm If $e=1$ then $\a_i=1$ if $R_{i+1}-R_i=1$ and $\a_i=(R_{i+1}-R_i)/2+1$
otherwise. 
\elm
\pf We have $R_{i+1}-R_i\geq -2e=-2$ and if $R_{i+1}-R_i$ is negative
then it is even. Thus $R_{i+1}-R_i$ is either $-2$ or it is $\geq
0$. If $R_{i+1}-R_i=-2e=-2$ or $R_{i+1}-R_i=2e-2=2-2e=0$ or if
$R_{i+1}-R_i\geq 2e=2$ then $\a_i=(R_{i+1}-R_i)/2+e=(R_{i+1}-R_i)/2+1$ by
Corollary 2.9(i). If $R_{i+1}-R_i=1$, which is odd and $<2e$, we have
$\a_i=R_{i+1}-R_i=1$ by Lemma 2.7(iii). \qed 

Since $\a_i$'s are uniquely defined by the $R_i$'s, condition (ii) of the main
theorem is superfluous since it follows from (i). Also,
$\od a_1\cdots a_i=\od b_1\cdots b_i$ so $\od a_1\cdots a_ib_1\cdots b_i$ is
even. So if $R_{i+1}-R_i\leq 1$ we have $d(a_1\cdots a_ib_1\cdots b_i)\geq
1\geq\a_i$. So condition (iii) is supefluous if $R_{i+1}-R_i\leq 1$. If
$R_{i+1}-R_i=2$ then $\a_2=2$, while if $R_{i+1}-R_i>2$ then $\a_i>2$. Thus in
these cases (iii) becomes $a_1\cdots a_ib_1\cdots b_i\in\fs\cup\D\fs$ if
$R_{i+1}-R_i=2$ and $a_1\cdots a_ib_1\cdots b_i\in\fs$ if
$R_{i+1}-R_i>2$. Finally it is easy to see that the condition
$\a_{i-1}+\a_i>2$ from 3.1(iv) is satisfied iff $R_{i-1}<R_{i+1}$ and the pair
$(R_i-R_{i-1},R_{i+1}-R_i)$ is different from $(0,1),(1,0),(1,1)$. So we have:

\btm Suppose that $F$ is 2-adic, $L\ap\[ a_1,\ldots,a_n\]$ and $K\ap\[
b_1,\ldots,b_n\]$ relative to good BONGs, $R_i=R_i(L)=\od a_i$,
$S_i=R_i(K)=\od b_i$ and $FL\ap FK$. Then $L\ap K$ if and only if the
following conditions hold: 

(i) $R_i=S_i$ for any $1\leq i\leq n$. 

(ii) For any $1\leq i\leq n-1$ we have
$a_1\cdots a_ib_1\cdots b_i\in\fs\cup\D\fs$ if $R_{i+1}-R_i=2$, and
$a_1\cdots a_ib_1\cdots b_i\in\fs$ if $R_{i+1}-R_i>2$. 

(iii) $[b_1,\ldots,b_{i-1}]\rep [a_1,\ldots,a_i]$ for any $1<i<n$
s.t. $R_{i-1}<R_{i+1}$ and $(R_i-R_{i-1},R_{i+1}-R_i)\neq
(0,1),(1,0),(1,1)$. 
\etm

\section{Remarks}

{\bf 1. The binary case}

If $L\ap\[\a,\b\]$ and $\eta\in\ooo$ then [B, 3.12] states that
$L\ap\[\eta\a,\eta\b\]$ iff $\eta\in g(a(L))=g(\h\b\a )$. 

The function $g:\aaa\z Sgp(\ooo/\ooo^2)$ was introduced in [B, Definition
6]. Here $Sgp\, H$ is the set of all subgroups of a group $H$. We recall the
definition of $g$.\footnote{In [B, Definition 6] there are some mistakes which
we corrected here.}

\no{\bf Definition} If $a=\pi^R\e\in\aaa$ and $d(-a)=d$ then: 

I If $R>2e$ then $g(a)=\ooo^2$.

II If $R\leq 2e$ then:
$$g(a)=\begin{cases}
\upo{R/2+e}&\text{ if }d>e-R/2\\ \upon{R+d}{-a}&\text{ if }d\leq e-R/2.
\end{cases}$$

The following lemma gives a more compact formula for $g(a)$.

\blm If $a\in\aaa$ and $\od a=R$ and $d(-a)=d$ then $g(a)=\upon{\a
(a)}{-a}$, where $\a (a)=\min\{ R/2+e,R+d\}$. 
\elm
\pf By [B, 3.16] we have $g(a)\sbq\N (-a)$. If $\eta\in\ooo$ then
$\eta\in g(a)$ iff $\eta\in\N (-a)$ and (I) If $R>2e$ then
$\eta\in\ooo^2$; (II) If $R\leq 2e$ then $d(\eta )\geq R+d$, if $d\leq
e-R/2$, and $d(\eta )\geq R/2+e$, if $d>e-R/2$. (See [B, Definition 6].) 

We have to prove that the conditions from (I) and (II) are equivalent
to $d(\eta )\geq\a (a)$. If $R>2e$ then $R+d>2e$ and $R/2+e>2e/2+e=2e$
so $\a (a)>2e$. Thus $d(\eta )\geq\a (a)$ is equivalent to
$\eta\in\ooo^2$. If $R\leq 2e$ then $d\leq e-R/2$ is equivalent to
$R+d\leq R/2+e$. Hence if $d\leq e-R/2$ then $\a (a)=R+d$ and if
$d>e-R/2$ then $\a (a)=R/2+e$. \qed

If $n=2$ then from [B, 3.12] we have $\[ a_1,a_2\]\ap\[\eta
a_1,\eta a_2\]$ iff $\eta\in g(a_2/a_1)$. By Lemma 5.1 this is
equivalent to $\eta\in\N (-a_1a_2)$ and $d(\eta )\geq\a
(a_2/a_1)$. The first condition is equivalent to the isometry of
quadratic spaces $[a_1,a_2]\ap [\eta a_1,\eta a_2]$, while the second
means $d(\eta )\geq\a (a_2/a_1)=\min\{
(R_2-R_1)/2+e,R_2-R_1+d(-a_1a_2)\} =\a_1(\[ a_1,a_2\] )$, which is
consistent with condition (iii) of the main theorem. 

\bff{\bf Remark} Since $\a (a_2/a_1)=\a_1(\[ a_1,a_2\] )$ we have by Lemma
5.1 $g(a_2/a_1)=\upon{\a_1(\[ a_1,a_2\] )}{-a_1a_2}$. Equivalently,
$g(a(L))=\upon{\a_1(L)}{-\det FL}$. 
\eff

\no{\bf 2. The formula for $\a_i$}

We will now show the heuristical method by which the invariants $\a_i$ were
found. We want to know, given that $L\ap\[ a_1,\ldots,a_n\]$ relative to a
good BONG and $1\leq i\leq n-1$, how much the product $a_1\cdots a_i$ can be
altered by a change of good BONGs. That is if $L\ap\[ b_1,\ldots,b_n\]$
relative to another good BONG we want to know how big the quadratic defect of
$(b_1\cdots b_i)/(a_1\cdots a_i)$ can be. So we are looking for a lower bound
$\a_i=\a_i(L)$ for $d(a_1\cdots a_ib_1\cdots b_i)$. 

For any $\eta\in g(a_{i+1}/a_i)$ we have $\[ a_i,a_{i+1}\]\ap\[\eta
a_i,\eta a_{i+1}\]$ so, by [B, Lemma 4.9(ii)], $L\ap\[
a_1,\ldots,a_{i-1},\eta a_i,\eta a_{i+1},a_{i+2},\ldots,a_n\]$. By this
change of BONGs $a_1\cdots a_i$ was changed by the factor $\eta$. We have
$\eta\in g(a_{i+1}/a_i)$ which, by Lemma 5.1, implies $d(\eta )\geq\a
(a_{i+1}/a_i)=\min\{
(R_{i+1}-R_i)/2+e,R_{i+1}-R_i+d(-a_ia_{i+1})\}$. (See Lemma 5.1.) This
lower bound can be further decreased if we decrease $d(-a_ia_{i+1})$. This can
be done by changing the good BONGs of $\[ a_1,\ldots,a_i\]$ and $\[
a_{i+1},\ldots,a_n\]$. If $\[ a_1,\ldots,a_i\]\ap\[ a'_1,\ldots,a'_i\]$ and $\[
a_{i+1},\ldots,a_n\]\ap\[ a'_{i+1},\ldots,a'_n\]$ then $d(-a_ia_{i+1})$ is
replaced by $d(-a'_ia'_{i+1})$. But $d(a_{i+1}a'_{i+1})\geq\a_1(\[
a_{i+1},\ldots,a_n\] )$. Also, by reason of determinant,
$a_1\cdots a_ia'_1\cdots a'_i\in\fs$ so
$d(a_ia'_i)=d(a_1\cdots a_{i-1}a'_1\cdots a'_{i-1})\geq\a_{i-1}(\[
a_1,\ldots,a_i\] )$. It follows that $d(-a'_ia'_{i+1})\geq\min\{
d(-a_ia_{i+1}),\a_{i-1}(\[ a_1,\ldots,a_i\] ),\a_1(\[ a_{i+1},\ldots,a_n\]
)\}$. Hence the new lower bound for $\eta$ is $\min\{
(R_{i+1}-R_i)/2+e, R_{i+1}-R_i+d(-a_ia_{i+1}), R_{i+1}-R_i+\a_{i-1}(\[
a_1,\ldots,a_i\] ),R_{i+1}-R_i+\a_1(\[ a_{i+1},\ldots,a_n\] )\}$. This leads
to the recursive formula $\a_i=\min\{
(R_{i+1}-R_i)/2+e, R_{i+1}-R_i+d(-a_ia_{i+1}), R_{i+1}-R_i+\a_{i-1}(\[
a_1,\ldots,a_i\] ),R_{i+1}-R_i+\a_1(\[ a_{i+1},\ldots,a_n\] )\}$ from
Corollary 2.5(ii). 

In the case $i=1$ and $n\geq 3$ the formula becomes $\a_1=\min\{
(R_2-R_1)/2+e,R_2-R_1+d(-a_1a_2),R_2-R_1+\a_1(\[ a_2,\ldots,a_n\]
)\}$. In the case $i=n-1$ and $n\geq 3$ we have $\a_{n-1}=\min\{
(R_n-R_{n-1})/2+e,R_n-R_{n-1}+d(-a_{n-1}a_n),R_n-R_{n-1}+\a_{n-2}(\[
a_1,\ldots,a_{n-1}\] )\}$. Finally if $i=1$ and $n=2$ then
$\a_1=\min\{ (R_2-R_1)/2+e,R_2-R_1+d(-a_1a_2)\}$. Starting with the
case $n=2$ it is easy to prove by induction that $\a_1=\min (\{
(R_2-R_1)/2+e\}\cup\{ R_{j+1}-R_1+d(-a_ja_{j+1})\mid 1\leq j<n\} )$ and
$\a_{n-1}=\min (\{ (R_n-R_{n-1})/2+e\}\cup\{
R_n-R_j+d(-a_ja_{j+1})\mid 1\leq j<n\} )$. By plugging $\a_{i-1}
(\[a_1,\ldots,a_i\] )=\min (\{ (R_i-R_{i-1})/2+e\}\cup\{
R_i-R_j+d(-a_ja_{j+1})\mid 1\leq j<i\} )$ and $\a_1(\[ a_{i+1},\ldots,a_n\]
)=\min (\{ (R_{i+2}-R_{i+1})/2+e\}\cup\{
R_{j+1}-R_{i+1}+d(-a_ja_{j+1})\mid i+1\leq j<n\} )$ in the recursive
formula for $\a_i$ we get the formula from Definition 1. (The extra
terms $R_{i+1}-R_i+(R_i-R_{i-1})/2+e$ and
$R_{i+1}-R_i+(R_{i+2}-R_{i+1})/2+e$ that appear are $\geq
(R_{i+1}-R_i)/2+e$ so they can be removed. )

Of course this is only a guess and does not constitute a proof. In
fact the relation $d(a_1\cdots a_ib_1\cdots b_i)\geq\a_i$ is only proved this
way in the particular case when $b_1,\ldots,b_n$ are obtained from
$a_1,\ldots,a_n$ through of succession of ``binary transformations'' of
the type $a_1,\ldots,a_n\z a_1,\ldots,\eta a_j,\eta a_{j+1},\ldots,a_n$ with
$1\leq j\leq n-1$ and $\eta\in g(a_{j+1}/a_j)$. It is not hard to prove that
conditions (i)-(iv) of the main theorem are necessary if $b_1,\ldots,b_n$ are
obtained this way. However, for the proof of the necessity in the general case
and for the proof of sufficiency the use of O'Meara's theorem is necessary. 

\no{\bf 3.} In the view of the previous remark there is the natural question
that asks whether, given that $L\ap\[ a_1,\ldots,a_n\]\ap\[ b_1,\ldots,b_n\]$
relative to good BONGs, there is always a succession of binary transformations
as defined above from $a_1,\ldots,a_n$ to $b_1,\ldots,b_n$. The answer to
this question is YES but only if we make the assumption that $F/\QQ_2$
is not totally ramified, i.e. that the residual field $\oo/\p$ has
more than 2 elements. 

If $|\oo/\p |=2$ we have the following counter-example. Let $0<d<2e$
be odd and let $R=2e-2d$ and $\e,\eta\in\ooo$ with $d(\e )=d$ and
$d(\eta )=2e-d$. It can be proved that $\[
1,-\pi^R\e,\e\eta,-\pi^R\eta\]\ap\[\eta,-\pi^R\e\eta,\e,-\pi^R\]$ but
one cannot go from $1,-\pi^R\e,\e\eta,-\pi^R\eta$ to
$\eta,-\pi^R\e\eta,\e,-\pi^R$ through binary transformations. 

E.g. if $F=\QQ_2$ and we take $d=1$, so $R=0$ and $\e =\eta =-1$, then
$\[ 1,1,1,1\]\ap\[ 7,7,7,7\]$. However from $1,1,1,1$ we can go through
binary transformations only to $a_1,a_2,a_3,a_4$, where an even number
of $a_i$'s belong to $\ooo^2$ and the rest to $5\ooo^2$. This happens
because $g(1)=g(5)=\ooo^2\cup 5\ooo^2$ so the only binary relations
involving $1$ and $5$ are $\[ 1,1\]\ap\[ 5,5\]$ and $\[ 1,5\]\ap\[
5,1\]$. Similarly from $7,7,7,7$ we can only go to $a_1,a_2,a_3,a_4$,
where an even number of $a_i$'s belong to $7\ooo^2$ and the rest to
$3\ooo^2$. 

\section*{References}

\hskip 6mm [B] C. N. Beli, Integral spinor norms over dyadic local
fields, J. Number Theory 102 (2003) 125-182.

[B1] C.N Beli, Representations of integral quadratic forms over dyadic local
fields, Electronic Research Announcements of the American Mathematical Society
12, 100-112, electronic only (2006).

[H] J. S. Hsia, Spinor norms of local integral rotations I, Pacific
J. Math., Vol. 57 (1975), 199 - 206.

[OM] O. T. O'Meara, Introduction to Quadratic Forms, Springer-Verlag,
Berlin (1963). 

\begin{center}
Institute of Mathematics of the Romanian Academy,

P.O. Box 1-764, RO-70700 Bucharest, Romania

email: raspopitu1@yahoo.com, Constantin.Beli@imar.ro
\end{center}

\end{document}